\documentclass[]{interact}

\usepackage{epstopdf}
\usepackage[caption=false]{subfig}
\usepackage{color}

\usepackage[numbers,sort&compress]{natbib}
\bibpunct[, ]{[}{]}{,}{n}{,}{,}

\theoremstyle{plain}
\newtheorem{theorem}{Theorem}[section]
\newtheorem{lemma}[theorem]{Lemma}
\newtheorem{corollary}[theorem]{Corollary}
\newtheorem{proposition}[theorem]{Proposition}

\theoremstyle{definition}
\newtheorem{definition}[theorem]{Definition}

\theoremstyle{remark}
\newtheorem{remark}{Remark}

\def\trt{^{\scriptscriptstyle T}}
\renewcommand{\int}{\mathop{\rm int}}

\usepackage[ruled,lined,boxed,norelsize]{algorithm2e}

\SetKwBlock{Repeat}{repeat}{end}

\begin{document}


\title{Diminishing Stepsize Methods for Nonconvex Composite Problems via Ghost Penalties: from the General to the Convex Regular Constrained Case
 }

\author{
\name{F. Facchinei\textsuperscript{a}\thanks{CONTACT F. Facchinei
Email: francisco.facchinei@uniroma1.it},  Vyacheskav Kungurtsev\textsuperscript{b},  
Lorenzo Lampariello\textsuperscript{c} 
and Gesualdo Scutari\textsuperscript{d} }
\affil{\textsuperscript{a}Sapienza University of Rome, Department of Computer, Control, and Management Engineering Antonio Ruberti, Rome, Italy; \textsuperscript{b}Czech Technical University in Prague, Department of Computer Science, Faculty of Electrical Engineering, Prague, Czech Republic; \textsuperscript{c}Roma Tre University, Department of Business Studies, Rome, Italy; \textsuperscript{d}Purdue University, School of Industrial Engineering,  West-Lafayette, IN, USA} }

\maketitle

\begin{abstract}
In this paper we first extend the diminishing stepsize method for nonconvex constrained problems presented in \cite{facchinei2017ghost} to deal with equality constraints and a nonsmooth objective function of composite type. We then consider the particular case in which the constraints are convex and satisfy a standard constraint qualification and show that in this setting the algorithm can be considerably simplified, reducing the computational burden of each iteration.
\end{abstract}

\begin{keywords}
Constrained optimization; Nonconvex optimization; Composite optimization, Diminishing stepsize.
\end{keywords}

\section{Introduction}\label{Sec:Intro}
We consider the  nonconvex constrained optimization problem
\begin{equation}\label{eq:startpro}
\begin{array}{cl}
\underset{x}{\mbox{minimize}} & f(x) + q(x) \\
\mbox{s.t.} & g(x) \le 0\\[5pt]
& h(x) = 0\\[5pt]
& x \in K,
\end{array}\tag{P}
\end{equation}
where $K \subseteq \mathbb R^n$ is a nonempty closed and convex set, and $f : \mathbb R^n \to \mathbb R,$  $g : \mathbb R^n \to \mathbb R^{m}$, $h : \mathbb R^n \to \mathbb R^p$ are C$^{1,1}$ (i.e., continuously differentiable with locally Lipschitz continuous gradients) functions on an open set containing $K,$ while $q : \mathbb R^n \to \mathbb R$ is convex on $K$ and locally Lipschitz continuous on an open set containing $K$.
In \cite{facchinei2017ghost}, building on an extended SQP-like approach, we analyze the first Diminishing Stepsize Method (DSM) for a general optimization problem with a nonconvex objective function and nonconvex constraints; we refer the reader to \cite{facchinei2017ghost} for a detailed discussion on DSMs in nonconvex settings.
However, the problem considered in \cite{facchinei2017ghost} includes neither equality constraints nor a possibly nonsmooth term $q(x)$ in the objective function. The main contributions of the present paper are
\begin{enumerate}
\item the study of the DSM for the  more general problem \eqref{eq:startpro} including equality constraints and a possibly nonsmooth term $q(x)$ in the objective;
\item the development of a new DSM for \eqref{eq:startpro} whenever constraints are \emph{convex} (with $f$ still possibly nonconvex), which results in a simpler algorithm than the one in (1) while maintaining global convergence guarantees.
\end{enumerate}
Regarding (1), the importance of the (possibly) nonsmooth term $q$ (for example the $\ell_1$ norm) cannot be overestimated: in the past ten years, objective functions including such terms have become pervasive, especially in the domains of machine learning and statistics. The presence of equality constraints of course also further widens the domain of applicability of our DSM. It should be remarked that equality constraints pose some subtle technical challenges making the extension of the results in \cite{facchinei2017ghost} to equality constrained problems much less immediate than might be anticipated.  Passing to the second contribution, we notice that, similarly to \cite{facchinei2017ghost}, the main cost of the new  DSM (cf. Section \ref{Sec:convergence}) is solving at each iteration two convex subproblems, whose  computational complexity  depends on the specific (inner) algorithm that is chosen to address them. However, if the constraints in~\eqref{eq:startpro} are known to be regular and convex, a simplified variant of the algorithm can be devised to good effect, which is more reminiscent of classical SQP methods and now only requires the solution of one strongly convex subproblem at each iteration; this algorithm is discussed in Section \ref{sec:newalgo}.

The algorithms we focus on are DSMs. At each iteration $\nu$,  a direction 
 $d(x^\nu)$ is computed  by solving a (surrogate) strongly convex optimization subproblem whose definition may require the solution of a further auxilary convex subproblem, actually a linear program if the appropriate choices are made. These subproblems are described  in Section \ref{Sec:dirfind}, where their properties are also analyzed in detail.
A step-size $\gamma ^\nu$ scales this direction so that the iterative sequence satisfies
\begin{equation}\label{eq:method}
x^{\nu+1} \,=\, x^\nu + \gamma^\nu d(x^\nu),
\end{equation}
with a stepsizes $\gamma^\nu$ such that the classical conditions
\begin{equation}\label{eq:gamma}
\lim_{\nu\to \infty} \gamma^\nu = 0\qquad \mbox{\rm and}\qquad \sum_{\nu=0}^\infty \gamma^\nu\, =\,
\infty
\end{equation}
hold. The convergence properties of this algorithm are analyzed in Section \ref{Sec:convergence}. Note that this algorithm does not require any assumptions, such as nonemptiness of the feasible set or constraint qualifications, to be well-defined and enjoy global convergence from arbitrary starting points to some limit point with desirable features. More specifically, the method is shown to be (subsequentially) convergent to points satisfying a (rather standard) generalized stationarity condition; the specific definition of generalized stationarity, along with its associated relevant properties, is discussed in Section \ref{sec:stationary}.
In the final part of the paper, see Section \ref{sec:newalgo}, we discuss the specific case in which the constraints are convex and a standard Mangasarian-Fromovitz-type constraint qualification is satisfed. For this setting, we introduce an alternative method, simpler than the one  discussed in Sections \ref{Sec:dirfind} and \ref{Sec:convergence}, and show its global (subsequential) convergence to KKT points.

We refer the reader to \cite{facchinei2017ghost} for a detailed literature review and comparison of the class of approaches we study with other algorithms; \cite{facchinei2017ghost} is also a good entry point for a discussion about generalized stationarity concepts.

We conclude this introduction by mentioning that an important  theoretical tool  we put forth  in our convergence analysis is the classical  penalty function
\[ W(x; \varepsilon) \triangleq f(x) + q(x) + \frac{1}{\varepsilon} \varphi(x),
 \]
where $\varepsilon$ is a positive penalty parameter, and
\begin{equation}\label{eq:feasmeas}
\varphi(x) \triangleq \max_{i,j}\{g_i(x)_+, \, |h_j(x)|\},
\end{equation}
where $\alpha_+ \triangleq  \max \{0,\alpha\}$.
The penalty function acts as a Lyapunov function in the convergence analysis in Section \ref{sec:newalgo}, when a constraint qualification holds. However, when analyzing the more general case considered in Section \ref{Sec:convergence}, the role of $W$ in the convergence analysis becomes more complex and is no longer that of a classical Lyapunov function. Note also that, in our analysis, the penalty function turns out to be just a theoretical tool and only enters in our convergence proof, and no penalty parameter needs to be actually computed in the algorithm itself. For these reasons, we refer to $W$ as a {\em ghost} penalty function.

\section{Generalized Stationarity}\label{sec:stationary}
In this section we introduce the concept of \emph{generalized stationarity}. As we study the convergence properties of our algorithm without assuming constraint qualifications or even nonemptiness of the feasible set 
of~\eqref{eq:startpro}, we need to define what sort of stationarity is possible and desirable to achieve in the case
that no limit point is KKT or even feasible. The concept of generalized stationarity, which characterizes stationarity across
the full taxonomy of possible desirable limit points, is discussed in more detail in~\cite{facchinei2017ghost} and the references
therein, and here we extend those considerations to account for the additional structural components of~\eqref{eq:startpro}.

Let us  denote the feasible set of \eqref{eq:startpro} by
\[
{\mathcal{X}} \triangleq \left\{x \in \mathbb R^n \, : \, g(x) \le 0, \, h(x) = 0, \, x \in K\right\}.
\]
The general constrained problem \eqref{eq:startpro} can be viewed as a combination of two problems: (i) the feasibility one, i.e., the problem of finding a feasible point; and (ii) the problem of finding a local minimum point of the objective function over the feasible set. Consistently, stationary solutions in a generalized sense are points that are either stationary for \eqref{eq:startpro} or for the following violation-of-the-constraints optimization problem:
\begin{equation}\label{eq:feaspro}
\begin{array}{cl}
\underset{x}{\mbox{minimize}} & \displaystyle \varphi(x),  \\
& x \in K
\end{array}
\end{equation}
where we recall that $\varphi(x)$, defined in \eqref{eq:feasmeas}, measures the degree of infeasibility. The KKT system for problem  \eqref{eq:startpro} is
\begin{equation}\label{eq:KKTstartpro}
\begin{array}{c}
0 \in \nabla f(x) + \partial q(x) + \nabla g(x) \xi + \nabla h(x) \pi + N_K(x) \\[5pt] 
0\leq \xi \perp g(x) \leq 0\\[5pt] 
h(x) =0\\[5pt] 
x\in K,
\end{array}
\end{equation}
where $\nabla f(x)$ is the gradient of $f$ in $x$, $\nabla g(x)$ and $\nabla h(x)$ are the transposed Jacobians of $g$ and $h$ evaluated at $x$, $\partial q$ is the subdifferential of $q$ at $x$, and $N_K(x)$ is the normal cone to $K$ at $x$. The vectors $\xi$ and $\pi$ appearing in~\eqref{eq:KKTstartpro}
are KKT \emph{multipliers} and an equivalent formulation of the condition is the existence of some pair $(\xi,\pi)$ in the set,
\begin{equation*}\label{eq:norm}
\begin{array}{rcl}
M_1(x) & \triangleq & \Big\{(\xi, \pi) \, \left| \right. \, \xi \in N_{\mathbb R^m_-}(g(x)), \, \pi \in N_0(h(x)),\\[5pt]
& & \hspace{45pt} 0 \in \nabla f(x) + \partial q(x) + \nabla g(x) \xi + \nabla h(x) \pi + N_K(x) \Big\}\\[5pt]
& = & \Big\{(\xi, \pi) \, \left| \right. \, \xi \in N_{\mathbb R^m_-}(g(x)), \, h(x) = 0, \\[5pt]
& & \hspace{45pt} 0 \in \nabla f(x) + \partial q(x) + \nabla g(x) \xi + \nabla h(x) \pi + N_K(x) \Big\}.
\end{array}
\end{equation*}
Note that indeed the KKT conditions are satisfied at a point $x$ if and only if $M_1(x) \neq \emptyset$. We also introduce the set of ``abnormal" multipliers
\begin{equation*}\label{eq:abn}
\begin{array}{rcl}
M_0(x) & \triangleq & \Big\{(\xi, \pi) \, \left| \right. \, \xi \in N_{\mathbb R^m_-}(g(x) - \varphi(x) e^m), \, \pi \in N_{\varphi(x) \mathbb B^p_\infty}(h(x)),\\[5pt]
 & & \hspace{45pt} 0 \in \nabla g(x) \xi + \nabla h(x) \pi + N_K(x) \Big\},
\end{array}
\end{equation*}
where $\mathbb B^p_\infty$ is the closed unit ball in $\mathbb R^p$ associated with the infinity-norm, $e^m\in \mathbb{R}^m$ is the vector of all ones, and, again,  $N_{\mathbb R^m_-}(y)$, $N_{\varphi(x) \mathbb B^p_\infty}(z)$ and $N_K(x)$ are the  normal cones to the convex sets $\mathbb R^m_-$, $\varphi(x) \mathbb B^p_\infty$ and $K$ at $y$, $z$ and $x$, respectively. 
Let  $\hat x$  be a local minimum point of \eqref{eq:startpro}, then it is well-known that either $M_1(\hat x) \neq \emptyset,$ (the point is a KKT point) or $M_0(\hat x) \neq \{0\}$ (the point is a Fritz-John point), or both. Otherwise, i.e., if $\hat x \in K$ is stationary but not feasible, in view of the regularity of the functions involved, then the appropriate stationarity condition is the one for problem \eqref{eq:feaspro}, i.e.,
\begin{equation}\label{eq:statfeaspro}
0 \in \partial \varphi(\hat x) + N_K(\hat x),
\end{equation}
which is equivalent to $M_0(\hat x) \neq \{0\}$. Hence, the (generalized) stationarity criteria for the original problem \eqref{eq:startpro} can naturally be specified by using the  sets $M_1$ and $M_0,$ as detailed in Definition \ref{df:stationarity}.
\begin{definition}\label{df:stationarity}
A point $\hat x \in K$ is, for problem \eqref{eq:startpro},
\begin{enumerate}
\item[$\bullet$] a KKT solution if $g(\hat x) \le 0,$ $h(\hat x) = 0$ and $M_1(\hat x) \neq \emptyset;$
\item[$\bullet$] a Fritz-John (FJ) solution if $g(\hat x) \le 0,$ $h(\hat x) = 0$ and $M_0(\hat x) \neq \{0\};$
\item[$\bullet$] an External Stationary (ES) solution if $g_i(\hat x) > 0$ or $h_j(\hat x) \neq 0$ for at least one $i \in \{1, \ldots, m\}$ or one $j \in \{1, \ldots, p\}$, and $M_0(\hat x) \neq \{0\}.$
\end{enumerate}
We call $\hat x \in K$ a stationary solution of \eqref{eq:startpro} if any of these cases occurs.
\end{definition}
The constraint qualification (CQ) we use is the Mangasarian-Fromovitz CQ, suitably extended to (possibly) infeasible points.
\begin{definition}\label{df:emfcq}
We say that the extended Mangasarian-Fromovitz Constraint Qualification (eMFCQ) holds at $\hat x \in K$ if
\begin{equation*}\label{eq:emfcq}
M_0(\hat x) = \{0\}.
\end{equation*}
\end{definition}
If $\hat x\in \mathcal{X}$   and $K = \mathbb R^n$, this condition reduces to the classical MFCQ and in turn, whenever the constraints are convex, it is well-known that the MFCQ is equivalent to Slater's CQ, i.e. to the existence of a point $\tilde x$ such that $g(\tilde x) < 0$ and $h(\tilde x) = 0$.
Below, we state a result that extends a standard property of the classical MFCQ for feasible points.

\begin{proposition}\label{th:cqlocal}
If the eMFCQ holds at $\hat x \in K$, then there exists a neighborhood $\mathcal V$ of $\hat x$ such that, for every $x \in K \cap \mathcal V$, the eMFCQ is satisfied.
\end{proposition}
\proof{}
If $\hat x \in K$ is feasible, this is a rather classical result: suffice it to reason by contradiction and to rely on the outer semicontinuity properties of the normal cone mappings $N_{\mathbb R^m_-}$, $N_K$ (see \cite[Proposition 6.6]{RockWets98}) and $N_{\varphi \mathbb B_\infty^p}$ (refer to point (iv) in the forthcoming Lemma \ref{th:prelresfeas}). If $\hat x \in K$ is not feasible, the condition $M_0(\hat x) = \{0\}$ implies that
$\hat x$ is not a stationary point for the feasibility problem \eqref{eq:feaspro}, i.e. $0\not \in \partial  \varphi(\hat x) + N_K(\hat x)$. The assertion then follows from the outer semicontinuity and local boundedness of the subdifferential mapping $\partial \varphi(\bullet)$ and by, again, the outer semicontinuity properties of the set valued mapping $N_K$ (see \cite{RockWets98} for the definition of outer semicontinuity). \hfill \endproof

\section{Algorithmic Scheme and Preliminary Results}\label{Sec:dirfind}
The key step in our algorithm is the computation of the direction $d(x^\nu)$ along which the update is performed. Specifically,   
at each iteration,  we move from the current iterate $x^\nu$ along the direction $d(x^\nu)$ with a stepsize $\gamma^\nu$ satisfying \eqref{eq:gamma}. We compute $d(x^\nu)$ as the solution of a strongly convex approximation of the original (possibly) nonconvex problem, with the approximating subproblem being reminiscent of (actually a generalization of) that of classical SQP methods. More precisely, given a point $x \in K$ (which will actually  be the current iterate $x^\nu$ in the algorithm), $d(x)$ is the unique solution of the following strongly convex optimization problem:
\begin{equation}\label{eq:p_k}
\begin{array}{cl}
\underset{d}{\mbox{minimize}} &  \tilde f(d;x) + q(x + d)\\
\mbox{s.t.} & \tilde g(d; x) \le \kappa(x) e^m\\[5pt]
& -\kappa(x) e^p \le \tilde h(d; x) \le \kappa(x) e^p,\\[5pt]
& \|d\|_\infty \le \beta,\\[5pt]
& d \in K - x
\end{array}\tag{P$_{x}$}
\end{equation}
where 
$\beta$ is a user-chosen positive constant, and the constraint $\|d\|_\infty \le \beta$ is introduced to prevent the direction $d(x^\nu)$ from becoming too large.  Moreover, $\tilde f$ is a strongly convex surrogate of function $f$, while $\tilde g_i$s are convex surrogate of the original constraint functions $g_i$ (see Assumption A below for the conditions these surrogates must satisfy), and
\[
\tilde h(d; x) \triangleq h(x) + \nabla h(x)\trt d.
\]
Finally,    the term
$\kappa(x)$ in the surrogate constraints is defined, for every $x \in K$, as  \begin{equation}\label{eq:cap}
\begin{array}{rcl}
\kappa(x) & \triangleq & (1 - \lambda) \max_{i,j}\{g_i(x)_+, \, |h_j(x)|\}\\[5pt]
& & + \lambda \min_d \left\{\max_{i, j}  \{\tilde g_i(d; x)_+, \, |\tilde h_j(d; x)| \} \, | \, \|d\|_\infty \le \rho, \, d \in K - x\right\},
\end{array}
\end{equation}
with $\lambda \in (0,1)$ and  $\rho \in (0, \beta)$. The goal of  $\kappa(x)$ is  making the feasible set of \eqref{eq:p_k} always nonempty.  
To computation of $\kappa(x)$ therefore requires one to calculate the optimal value of the convex problem
\begin{equation}\label{eq:capoptpro}
\min_d \left\{\max_{i, j}  \{\tilde g_i(d; x)_+, \, |\tilde h_j(d; x)| \} \, | \, \|d\|_\infty \le \rho, \, d \in K - x\right\}\end{equation}
that is always solvable because its feasible set is nonempty and compact. 
Note that in the most common case in which linear approximations are used  for the inequality constraints, this problem can easily be reformulated as an LP and hence efficiently solved; nevertheless, the computation of 
$\kappa(x)$ is a somewhat expensive task, the more so when nonlinear approximations are used for $g$. In the last section of this paper we see that under some additional assumptions the burden of this computation can be avoided altogether.
Furthermore, note that if $x$ is feasible for \eqref{eq:startpro}, then $\kappa(x) = 0$. 

 The method we propose is summarized in Algorithm \ref{algoBasic} below. 
\begin{algorithm}[h!]
\KwData{$\gamma^\nu \in (0,1]$ such that \eqref{eq:gamma} holds, $x^{0} \in K$,  $\nu \longleftarrow 0$\;}
\Repeat{
{\nlset{(S.1)} \If{$x^{\nu}$ {\emph{is a generalized stationary point of }} \eqref{eq:startpro}}{{\bf{stop}} and {\bf{return}} $x^\nu$\;} \label{S.11}}
{\nlset{(S.2)} compute $\kappa(x^\nu)$ and the solution $d(x^\nu)$ of problem (P$_{x^\nu}$)\;} \label{S.12}
{\nlset{(S.3)} set $x^{\nu+1}=x^{\nu}+\gamma^{\nu}d(x^\nu)$, $\nu\longleftarrow\nu+1$\;} \label{S.13}}{}
\caption{\label{algoBasic} DSM Algorithm for \eqref{eq:startpro}}
\end{algorithm}
A few remarks are in order. Subproblem \eqref{eq:p_k} is a generalization, along the lines explored in \cite{facchinei2017feasible,scutari2014decomposition}, of the direction finding subproblem considered in
\cite{burke1989sequential}, to which it reduces when the classical quadratic/linear approximations are used for 
$\tilde f$ and $\tilde g$:
\begin{equation}\label{eq:quad_lin approx}
\tilde f (d; x) \triangleq \nabla f(x)\trt d + \frac{1}{2} \| d\|_B^2; \qquad \quad
\tilde g (d; x) \triangleq g(x) + \nabla g(x)\trt d.
\end{equation}
where $B$ is some positive definite symmetric matrix. 
Note that if these approximations are employed and we set $\kappa(x) =0$ and $\beta = +\infty$, \eqref{eq:p_k} boils down to the classical SQP-type subproblem. In Section \ref{sec:newalgo} we shall see that under certain conditions,
setting $\kappa(x)=0$ and $\beta=+\infty$ is appropriate.  For the time being we adopt the approach in
\cite{burke1989sequential} by taking $\kappa(x)$ not necessarily zero and $\beta < +\infty$ in order to guarantee the existence and continuity of the solution mapping $d(x)$. In addition, we introduce the use of general 
approximations $\tilde f$ and $\tilde g$: this may be very convenient in practice by allowing flexibility in 
tailoring the direction finding subproblem to the problem at hand and to exploit any available specific structure 
in \eqref{eq:startpro} amenable to fast computation. It is clear that we implicitly suppose that the solution of subproblem \eqref{eq:p_k} is 
``easy" and, in any case, simpler than the original problem \eqref{eq:startpro}. We do not insist on this
point because it is very dependent on the choice of $\tilde f$ and $\tilde g$ which, in turn,  is guided by the
original problem (P). It is worth mentioning that the use of models that go beyond the standard
quadratic/linear one in  optimization is steadily emerging in the literature, motivated, on the one hand, by
the advances in the efficient solution of more complex subproblems than the classical quadratic ones and, on the 
other hand, by the desire of faster convergence rates, see for example the discussion in Section 3 of 
\cite{martinez2017high}.


In the sequel we denote by $\widetilde {\mathcal{X}}(x)$ and $d(x)$ the convex feasible set and the unique solution of subproblem \eqref{eq:p_k}, respectively, i.e.
\[
\begin{array}{rcl}
\widetilde {\mathcal{X}}(x) & \triangleq & \Big\{d \in \mathbb R^n \, : \, \tilde g(d; x) \le \kappa(x) \, e^m, \, -\kappa(x) e^p \le \tilde h(d; x) \le \kappa(x) e^p,\\[5pt]
& & \hspace{50pt} \|d\|_\infty \le \beta, \, d \in K - x \Big\},\\[10pt]
d(x) & \triangleq & \arg\min_d \{\tilde f(d; x) \, | \, d \in \widetilde{\mathcal{X}}(x)\},
\end{array}
\]
and we equivalently write the constraints $-\kappa(x) e^p \le \tilde h(d; x) \le \kappa(x) e^p$ and  $\|d\|_\infty \le \beta$ as $\tilde h(d; x) \in \kappa(x) \mathbb B^p_\infty$ and $d \in \beta \mathbb B^n_\infty$, respectively. 

For our approach to be legitimate and lead to useful convergence results, we obviously need to make assumptions on the surrogate functions $\tilde f$ and $\tilde g$.

\medskip \noindent {\bf{Assumption A}}

\medskip

\noindent {\textit{Let $O_d$ and $O_x$ be open neighborhoods of $\beta \mathbb B^n_\infty$ and $K$, respectively, and $\tilde{f}:O_d\times O_x\rightarrow\mathbb{R}$ and $\tilde g_i:\mathbb R^n\times O_x\rightarrow\mathbb R$, for every $i = 1, \ldots, m$, be continuously differentiable on $O_d$ with respect to the first argument and satisfy}}
\begin{description}
\item [{\textit{A1)}}] {\textit{$\tilde{f}(\bullet;x)$ is a strongly convex function on $O_d$ for every $x \in K$, with modulus of strong convexity $c >0$ independent of $x$;}}

\item [{\textit{A2)}}] {\textit{$\tilde f(\bullet;\bullet)$ is continuous on $O_d \times O_x$;}}

\item [{\textit{A3)}}] {\textit{$\nabla_1 \tilde{f}(\bullet;\bullet)$ is continuous $O_d \times O_x$;}}

\item [{\textit{A4)}}] {\textit{$\nabla_1 \tilde f(0; x) = \nabla f(x)$ for every $x \in K$;}}

\item [{\textit{A5)}}] {\textit{$\tilde g_i(\bullet;x)$ is a convex function on $O_d$ for every $x \in K$;}}

\item [{\textit{A6)}}] {\textit{$\tilde g_i(\bullet;\bullet)$ is continuous on $\mathbb R^n \times O_x$;}}

\item [{\textit{A7)}}] {\textit{$\tilde g_i(0;x) = g_i(x)$ for every $x \in K$;}}

\item [{\textit{A8)}}] {\textit{$\nabla_1 \tilde g_i(\bullet;\bullet)$ is continuous on $O_d \times O_x$;}}

\item [{\textit{A9)}}] {\textit{$\nabla_1 \tilde g_i(0;x)=\nabla g_i(x),$ for every $x \in K$;}}
\end{description}
{\textit{where $\nabla_1 \tilde f(u;x)$ and $\nabla_1 \tilde g_i(u;x)$ denote the partial gradient of $\tilde f(\bullet;x)$ and $\tilde g_i(\bullet;x)$ evaluated at $u$.}}

These conditions are  certainly satisfied if we use the classical surrogates \eqref{eq:quad_lin approx}, but they allow us to cover a much wider array of approximations, both for $f$ and for $g$; we refer the reader to \cite{facchinei2017feasible,scutari2014decomposition} as good sources of examples of nonlinear surrogates $\tilde f$ and $\tilde g$ satisfying Assumption A.
Note that, under Assumption A, 
Algorithm \ref{algoBasic} is always well defined and $d(x^\nu)$ exists and is unique.

\subsection{Main Properties of Subproblem \eqref{eq:p_k}}\label{sec:subpro} In this section, we state the main properties of $\kappa (x)$  and problem \eqref{eq:p_k}. 

Since $\kappa (x)$ is always nonnegative, it restores feasibility by enlarging the range of admissible constraint function values.
Indeed, the point $\tilde d$ where a minimum is reached in the optimization problem in  \eqref{eq:cap} is easily seen to be always  feasible for \eqref{eq:p_k}.
Moreover, for every $x\in K$, the following relations hold thanks to Assumption A:
{\small{\begin{equation}\label{eq:kbounds1}
	\min_d \left\{\max_{i, j}  \{\tilde g_i(d; x)_+, \, |\tilde h_j(d; x)| \} \, | \, \|d\|_\infty \le \rho, \, d \in K - x\right\} \le \kappa(x) \le \max_{i,j}\{g_i(x)_+, \, |h_j(x)|\},  
\end{equation}}}
and   
\begin{equation}\label{eq:kbounds2}
\begin{array}{c}
\kappa(x) = \max_{i,j}\{g_i(x)_+, \, |h_j(x)|\}\\[5pt]
\Updownarrow\\[5pt]
	\max_{i,j} \{g_i(x)_+, |h_j(x)|\} = \min_d \left\{\max_{i, j}  \{\tilde g_i(d; x)_+, \, |\tilde h_j(d; x)| \} \, | \, \|d\|_\infty \le \rho, \, d \in K - x\right\}\\[5pt]
	\Updownarrow\\[5pt]
	0 \in \arg\min_d \left\{\max_{i, j}  \{\tilde g_i(d; x)_+, \, |\tilde h_j(d; x)| \} \, | \, \|d\|_\infty \le \rho, \, d \in K - x\right\}.
\end{array}   
\end{equation}
Also, $\kappa(x) = 0$ if and only if $0 = \max_{i,j}\{g_i(x)_+, |h_j(x)|\} = \min_d \{\max_{i,j} \{\tilde g_i(d;x)_+, |\tilde |h_j(d;x)|\} | \|d\|_\infty \le \rho, d \in K - x\}$.

In the remaining part of this section we give some highly technical results needed in the convergence analysis. The developments  are along lines similar to those considered in \cite[Section 3.2]{facchinei2017ghost} with, however, the additional hurdle of the equality constraints and the nonsmooth term $q(x)$ that need a specific, in some cases non trivial, treatment. 

In Lemma \ref{th:prelresfeas}, we establish some preliminary properties concerning the feasible set of problem \eqref{eq:p_k}.
\begin{lemma}\label{th:prelresfeas}
\begin{enumerate}
\item[(i)]
For every $\hat x \in K$, and for every $\alpha > 0$ and $d \in \alpha \mathbb B_\infty^n \cap (K - \hat x)$, the constraint qualification
\begin{equation}\label{eq:furcq}
[-N_{\alpha \mathbb B^n_\infty}(d)] \cap N_{K-\hat x}(d) = \{0\}
\end{equation}
holds and, in turn, $N_{\alpha \mathbb B^n_\infty \cap (K - \hat x)} (d) = N_{\alpha \mathbb B_\infty^n}(d) + N_{K - \hat x}(d)$;
\item[(ii)]
for every $\alpha > 0$, the set-valued mapping $\alpha \mathbb B_\infty^n \cap (K - \bullet)$ is continuous on  $K$ relative to $K$;
\item[(iii)]
letting $C \triangleq \{(d, x) \in \beta \mathbb B_\infty^n \times K \, : \, d + x \in K \}$, the set-valued mapping $N_{\beta \mathbb B_\infty^n \cap (K - \bullet)}(\bullet)$ is outer semicontinuous on $C$ relative to $C$;
\item[(iv)]
letting $\psi: \, \mathbb R^n \to \mathbb R$ be any nonnegative function that is continuous on $K$ relative to $K$, and $D \triangleq \{(u, x) \in \psi(x) \mathbb B_\infty^p \times K\}$, the set-valued mapping $N_{\psi(\bullet) \mathbb B_\infty^p}(\bullet)$ is outer semicontinuous on $D$ relative to $D$.  
\end{enumerate}
\end{lemma}
\proof{}
(i) Let $0 \neq \eta \in [- N_{\alpha \mathbb B_\infty^n}(d)] \cap N_{K- \hat x}(d)$. Thanks to the convexity of the sets $\alpha \mathbb B_\infty^n$ and $K - \hat x$, we have $- \eta\trt (v - d) \le 0$ $\forall v \in \alpha \mathbb B^n_\infty$ and $\eta\trt (y - d) \le 0$ $\forall y \in (K - \hat x)$. Choosing $y = 0 \in (K - \hat x)$, one gets the following contradiction:
\[
0 < \alpha \, \max_v \{-\eta\trt v \, | \, \|v\|_\infty = 1\} = -\eta\trt d \le 0,
\]
thus proving relation \eqref{eq:furcq}. As a consequence, the other claim in (i) follows from \cite[Theorem 6.42]{RockWets98}.

(ii)  The property is due to the continuity (relative to $K$) of the set-valued mapping $K - \bullet$ at every $x \in K$ and to the fact that $\alpha \mathbb B_\infty^n \cap (K - x) \neq \emptyset$ for every $x \in K$.

(iii) Suppose by contradiction that $(d^\nu, x^\nu) \underset{C}{\to} (\bar d, \bar x)$, $\ \eta^\nu \in N_{\beta \mathbb B_\infty^n \cap (K - x^\nu)}(d^\nu)$, $\eta^\nu \to \bar \eta$ with $\bar \eta \notin N_{\beta \mathbb B_\infty^n \cap (K - \bar x)}(\bar d)$. Hence, $\bar z \in \beta \mathbb B_\infty^n \cap (K - \bar x)$ exists such that $\bar \eta\trt (\bar z - \bar d) > 0$. By the inner semicontinuity relative to $K$ (see \cite[Chapter 5, Section B]{RockWets98} for the definition of inner semicontinuity) of $\beta \mathbb B_\infty^n \cap (K - \bullet)$ at $\bar x$, $z^\nu$ exists such that $z^\nu \to \bar z$ and $z^\nu \in \beta \mathbb B_\infty^n \cap (K - x^\nu)$. In turn, eventually we get $(\eta^\nu)\trt (z^\nu - d^\nu) > 0$ in contradiction to the inclusion $\ \eta^\nu \in N_{\beta \mathbb B_\infty^n \cap (K - x^\nu)}(d^\nu)$.

(iv) Taking into account the continuity of the set-valued mapping $\psi(\bullet) \mathbb B_\infty^p$, the proof follows the same line of reasoning as in (iii). 
\hfill  \endproof
The function $\kappa(x)$ is obviously continuous and, under the local Lipschitz continuity of $\tilde g(\bullet; \bullet)$ (which is  part of Assumption C to be introduced shortly), also locally Lipschitz continuous. This result has been shown in \cite{burke1989sequential} when
$\tilde g$ is a linear approximation and readily generalizes to the case of the surrogate $\tilde g$ we consider here.\begin{proposition}\label{eq:kappa lip}
Under Assumption A, $\kappa(\bullet)$ is  continuous on $K$ relative to $K$. If, in addition, $\tilde g(\bullet; \bullet)$ is locally Lipschitz continuous on $O_d \times O_x$, then $\kappa(\bullet)$ is also locally Lipschitz continuous on an open neighborhood of $K$.
\end{proposition}
\proof{}
The continuity of  $\kappa(\bullet)$ follows readily from the continuity (relative to $K$) of the set-valued mapping $\rho \mathbb B_\infty^n \cap (K - \bullet)$ at every $x \in K$: this in turn follows from (ii) in Lemma \ref{th:prelresfeas} with $\alpha = \rho$.

The Lipschitz continuity under the additional condition derives from, e.g., \cite[Theorem 3.1]{Rock85}. Suffice it to observe that the constraint qualification \eqref{eq:furcq} with $\alpha = \rho$ holds for every $x \in K$ and $d \in \rho \mathbb B_\infty^n \cap K - x$, and the problem \eqref{eq:capoptpro} in the definition of $\kappa$ is solvable for every $x$ in an open neighborhood of $K$. The latter claim is due to $\rho \mathbb B_\infty^n \cap (K - x) \neq \emptyset$ for every $x$ in an open neighborhood of $K$, since, for every $x \in K$, $0 \in \mathrm{int} (\rho \mathbb B_\infty^n) \cap (K - x)$, and in view of the continuity of the set-valued mapping $K - \bullet$.
\hfill  \endproof

The following technical lemma is very useful for the subsequent developments.
\begin{lemma}\label{th:prellem}
Under Assumption A, the following results hold for any $\hat x \in K$:
\begin{enumerate}
\item[(i)]
if $\max_{i,j}\{g_i(\hat x)_+, \, |h_j(\hat x)|\} > 0$ and $\kappa(\hat x) < \max_{i,j}\{g_i(\hat x)_+, \, |h_j(\hat x)|\}$, then, for all $\rho \in (0, \beta)$, there exists $d \in \mathrm{int} (\beta \mathbb B^n_\infty) \cap \mathrm{rel \, int}(K - \hat x)$ such that $\tilde g(d; \hat x) < \kappa(\hat x) e^m$ and $\kappa(\hat x) e^p < \tilde h(d; \hat x) < \kappa(\hat x) e^p$;

\item[(ii)]
if $\max_{i,j} \{g_i(\hat x)_+, |h_j(\hat x)|\} > 0$ and $\kappa(\hat x) = \max_{i,j} \{g_i(\hat x)_+, |h_j(\hat x)|\}$, then $\hat x$ is an ES point for \eqref{eq:startpro};

\item[(iii)]
if $\max_{i,j} \{g_i(\hat x)_+, |h_i(\hat x)|\} = 0$, then either $\hat x$ is a FJ point for \eqref{eq:startpro} or, for all $\rho \in (0, \beta)$, there exists $d \in \mathrm{int} (\beta \mathbb B^n_\infty) \cap \mathrm{rel \, int}(K - \hat x)$ such that $\tilde g(d; \hat x) < 0$, $\tilde h(d; \hat x) = 0$ and $\{0\} = \Big\{\pi \, | \, (\nabla h(\hat x) \pi)\trt w = 0, \forall w \in T_{\beta \mathbb B_\infty^n \cap (K - \hat x)}(v) \Big\}$ for every $v \in \widetilde{\mathcal X}(\hat x)$, where $T_{\beta \mathbb B_\infty^n \cap (K - \hat x)}(v)$ denotes the tangent cone to  ${\beta \mathbb B_\infty^n \cap (K - \hat x)}$ at $v$.
\end{enumerate}
\end{lemma}
\proof{}
(i) First, we recall that in this case, in view of the preliminary relations discussed at the beginning of Section \ref{sec:subpro}, $\kappa(\hat x) > 0$. Choosing $\hat d \in \arg\min_d \{\max_{i, j}  \{\tilde g_i(d; x)_+, \, |\tilde h_j(d; x)| \} \, | \, \|d\|_\infty \le \rho, \, d \in K - x\}$ with $\rho \in (0, \beta)$, we can infer $\tilde g_i(\hat d; \hat x) \le \min_d \{\max_{i,j}\{\tilde g_i(d; \hat x)_+, |h_j(d; \hat x)|\} \, | \, \|d\|_\infty \le \rho, d \in K - \hat x \}$ for every $i$, and $|\tilde h_j(\hat d; \hat x)| \le \min_d \left\{\max_{i,j}\{\tilde g_i(d; \hat x)_+, |h_j(d; \hat x)|\} \, | \, \|d\|_\infty \le \rho, d \in K - \hat x \right\}$ for every $j$, with $\hat d \in \rho \mathbb B_\infty^n \cap (K - \hat x)$. The claim follows in view of \eqref{eq:kbounds1} and \eqref{eq:kbounds2}, and by continuity since $\rho < \beta$.

(ii) By \eqref{eq:kbounds2}, equality $\kappa(\hat x) = \max_{i,j} \{g_i(\hat x)_+, |h_j(\hat x)|\}$ holds if and only if $d = 0$ solves the minimization problem in the definition of $\kappa$ and, in turn, $M_0(\hat x) \neq \{0\}$ since condition \eqref{eq:statfeaspro} holds at $\hat x$ by \eqref{eq:furcq} with $\alpha = \rho$, A7 and A9.

(iii) With $\max_{i,j} \{g_i(\hat x)_+, |h_j(\hat x)|\}$ being equal to zero, we have $\kappa(\hat x) = 0$, and $g(\hat x) \le 0$ and $h(\hat x) = 0$. If $M_0(\hat x) \neq \{0\}$, then, by definition, $\hat x$ is a FJ point for \eqref{eq:startpro} and the result holds.

Thus, let us suppose $M_0(\hat x) = \{0\}$. Following the same line of reasoning as in \cite[Exercise 6.39]{RockWets98}, in view of the regularity of the involved sets, we see preliminarily that this condition holds at $\hat x \in K$ if and only if
\begin{gather}
	\{0\} = \Big\{\pi \, | \, (\nabla h(\hat x) \pi)\trt w = 0, \forall w \in T_K(\hat x) \Big\}, \nonumber \\[7pt]
 \exists \, \hat d \,\in\, \mathrm{rel \, int} \, T_K(\hat x):\; \nabla g_i(\hat x)\trt \hat d < 0, \quad \forall i: g_i(\hat x) = 0, \qquad \nabla h(\hat x)\trt \hat d = 0. \label{eq:EMFCQ1}
\end{gather}
First, we show that condition
\begin{equation}\label{eq:linearindenh}
\{0\} = \Big\{\pi \, | \, (\nabla h(\hat x) \pi)\trt w = 0, \forall w \in T_{\beta \mathbb B_\infty^n \cap (K - \hat x)}(d) \Big\}	
\end{equation}
holds for any $d \in {\widetilde {\mathcal{X}}}(\hat x)$. If this were not the case, there would exist some $\tilde d \in {\widetilde {\mathcal{X}}}(\hat x)$ and $\pi \neq 0$ such that $(\nabla h(\hat x) \pi)\trt w = 0, \forall w \in T_{\beta \mathbb B_\infty^n \cap (K - \hat x)}(\tilde d)$. Thus, with $- \nabla h(\hat x)\pi \in N_{\beta \mathbb B_\infty^n \cap (K - \hat x)}(\tilde d)$, for some $\eta \in N_{\beta \mathbb B_\infty^n}(\tilde d)$ and $\zeta \in N_{K - \hat x}(\tilde d)$ it would hold that,
\begin{equation*}\label{eq:interlinearindenh}
0 = (\nabla h(\hat x) \pi)\trt (-\tilde d) + \eta\trt (-\tilde d) + \zeta\trt (-\tilde d)= \eta\trt (-\tilde d) + \zeta\trt (-\tilde d),
\end{equation*}
thanks to \cite[Theorem 6.42]{RockWets98} and observing that $- \tilde d \in T_{\beta \mathbb B_\infty^n \cap (K - \hat x)}(\tilde d)$ since $0 \in \beta \mathbb B_\infty^n \cap (K - \hat x)$. In turn,
\[
0 = \eta\trt \tilde d + \zeta\trt \tilde d \ge \rho \, \max_v \{\eta\trt v \, | \, \|v\|_\infty=1\} \ge 0,
\]
where the first inequality follows from $\zeta\trt \tilde d \ge 0$, in view of $0 \in \beta \mathbb B_\infty^n$. As a consequence, $\eta = 0$ and $\zeta \trt \tilde d = 0$, entailing $\zeta \in N_K(\hat x)$. Therefore, we would get
\[
0 = \nabla h(\hat x) \pi + \zeta,
\]
with $\pi \neq 0$ and $\zeta \in N_K(\hat x)$, in contradiction with the assumed condition $M_0(\hat x) = \{0\}$. Hence, \eqref{eq:linearindenh} is verified at any $d \in \widetilde {\mathcal{X}}(\hat x)$.    

Then, for those $j \in \{1, \ldots, m\}$ such that $g_j(\hat x) < 0$, we have $\tilde g_j(0; \hat x) = g_j(\hat x) < 0$; as for indices $k \in \{1, \ldots, m\}$ with $g_k(\hat x) = 0$, by \eqref{eq:EMFCQ1}, there exists $\hat d \in \mathrm{rel \, int} \, T_K(\hat x)$ such that
$$
0 > \nabla g_k(\hat x)\trt \hat d = \nabla_1 \tilde g_k(0; \hat x)\trt \hat d =  \lim_{\tau \downarrow 0} \frac{\tilde g_k(\tau \hat d; \hat x) - \tilde g_k(0; \hat x)}{\tau},
$$
as well as $\nabla h(\hat x)\trt \hat d = 0$. Taking the cue from the proof of \cite[Theorem 6.9]{RockWets98}, one can observe that $\mathrm{rel \, int} \, T_K(\hat x) = \left\{d \in \mathbb R^n \, | \, \exists \, \alpha > 0 \enspace \text{with} \enspace \hat x + \alpha d \in \mathrm{rel \, int} K \right\}$ due to \cite[Proposition 2.40]{RockWets98}. Hence, in view of \cite[Theorem 6.1]{rockafellar1970convex}, for every $\tau > 0$ sufficiently small, $\hat x + \tau \hat d \in \mathrm{rel \, int} K$ as well. The claim follows by continuity, observing that $\tilde g_i(\tau \hat d; \hat x) < 0$ for every $i$ and for any $\tau$ sufficiently small.
\hfill  \endproof
The quantity
\begin{equation}\label{eq:theta}
\begin{array}{r}
\theta(x)  \triangleq  \max_{i,j} \{g_i(\hat x)_+, |h_j(\hat x)|\} - \kappa(x) = \lambda \Big(\max_{i,j} \{g_i(\hat x)_+, |h_j(\hat x)|\}\\[5pt]
-\min_d \left\{\max_{i, j}  \{\tilde g_i(d; x)_+, \, |\tilde h_j(d; x)| \} \, | \, \|d\|_\infty \le \rho, \, d \in K - x\right\}\Big),
\end{array}
\end{equation}
with $\lambda\in (0,1)$, plays a key role in the previous lemma and in all the subsequent developments. As shown in the following proposition, $\theta$ turns out to be a stationarity measure for the violation-of-the-constraints problem \eqref{eq:feaspro}.
\begin{proposition}\label{th:theta}
Under Assumption A,
\begin{itemize}
\item[(i)]	
the nonnegative function $\theta(\bullet)$ is continuous on $K$ relative to $K$;
\item[(ii)]
$\theta(\hat x) = 0$ if and only if $\hat x$ is a stationary point for problem \eqref{eq:feaspro};
\item[(iii)]
we have, for every $x \in K$,
\begin{equation}\label{eq:thetadelta}
\theta(x) \le \left\|\begin{pmatrix}
\nabla g(x^\nu)\trt\\
\nabla h(x^\nu)\trt	
\end{pmatrix}\right\|_\infty \|d(x)\|.
\end{equation}
\end{itemize}
\end{proposition}
\proof{}
(i) The nonnegativity of $\theta$ follows readily from \eqref{eq:kbounds1} while 
continuity follows from Proposition \ref{eq:kappa lip}.

(ii) At any feasible point $\hat x$ of \eqref{eq:startpro}, $\theta(\hat x) = 0$; of course, every feasible point of \eqref{eq:startpro} is stationary for problem \eqref{eq:feaspro}. Consider now an infeasible point $\hat x$ for \eqref{eq:startpro} and suppose that $\theta(\hat x)=0$. By (ii) in Lemma \ref{th:prellem}, $\hat x$ turns out to be an ES point for \eqref{eq:startpro}. Hence, we are left to show that if $\hat x$ is an ES point for \eqref{eq:startpro}, then $\theta(\hat x) = 0$. For $\hat x$ to be ES, it is necessary and sufficient (see condition \eqref{eq:statfeaspro}) to have $M_0(\hat x) \neq \{0\}$ which in turn, by the Motzkin's alternative theorem (see e.g. \cite[2.5.2]{craven2012mathematical}), holds if and only if
\begin{equation}\label{eq:Motzk}
\begin{array}{c}
\nexists \, d \in T_K(\hat x) :  \nabla g_i(\hat x)\trt d < 0, \, \forall i \in I_+(\hat x) \triangleq \{i \, : \, g_i(\hat x) = \max_{i,j} \{g_i(\hat x)_+, |h_j(\hat x)|\}\},\\[5pt]
 \nabla h_j(\hat x)\trt d < 0, \, \forall j \in J_+(\hat x) \triangleq \{j \, : \, h_j(\hat x) = \max_{i,j} \{g_i(\hat x)_+, |h_j(\hat x)|\}\},\text{ and,}\\[5pt]
 \nabla h_j(\hat x)\trt d > 0, \, \forall j \in J_-(\hat x) \triangleq \{j \, : \, h_j(\hat x) = -\max_{i,j} \{g_i(\hat x)_+, |h_j(\hat x)|\}\}.	
 \end{array}
\end{equation}
Suppose by contradiction that $\theta(\hat x) > 0$. Then, noting that $d\in K-\hat x$ implies $d\in T_K(\hat x)$, Lemma \ref{th:prellem} (i) states that $d \in T_K(\hat x)$ exists such that
$\tilde g_i(d;\hat x) < \kappa(\hat x)$ for all $i \in I_+(\hat x)$ and $|\tilde h_j(d;\hat x)|<\kappa(\hat x)$ for all $j \in J_+(\hat x) \cup J_-(\hat x)$. But then, using A5, A7, and A9, we can write, for any $d$,
\[
\begin{array}{c}
\max_{i,j} \{g_i(\hat x)_+, |h_j(\hat x)|\}\}  \,>\, \kappa(\hat x) \,> \,\tilde g_i(d;\hat x)  \,\geq\, \tilde g_i(0;\hat x) + \nabla_1 \tilde g_i(0;\hat x)\trt(d-0)\\[5pt]
\geq\,  g_i(\hat x) + \nabla g_i(\hat x)\trt d,\\[10pt]
\max_{i,j} \{g_i(\hat x)_+, |h_j(\hat x)|\}\}  \,>\, \kappa(\hat x) \,> \,\tilde h_j(d;\hat x)  =  h_j(\hat x) + \nabla h_j(\hat x)\trt d,\\[10pt]
-\max_{i,j} \{g_i(\hat x)_+, |h_j(\hat x)|\}\}  \,<\, - \kappa(\hat x) \,< \,\tilde h_j(d;\hat x)  =  h_j(\hat x) + \nabla h_j(\hat x)\trt d,
\end{array}
\]
for every $i \in I_+(\hat x)$, $j \in J_+(\hat x)$ and $j \in J_-(\hat x)$, respectively. In turn, we get a contradiction to \eqref{eq:Motzk}.

(iii) Furthermore,
\begin{equation*}\label{eq:phicappa}
\begin{array}{rcl}
0 & \le & \theta(x^\nu) = \displaystyle \max_{i,j} \{g_i(x^\nu)_+, |h_j(x^\nu)|\}\} - \kappa(x^\nu)\\[5pt]
& \overset{(a)}{\le} & \displaystyle \max_{i,j} \{g_i(x^\nu)_+, |h_j(x^\nu)|\}\} - \max_{i,j}\{\tilde g_i(d(x^\nu); x^\nu)_+, |\tilde h_j(d(x^\nu); x^\nu)|\}\\[5pt]
& \overset{(b)}{\le} & \displaystyle \max_{i,j} \{g_i(x^\nu)_+, |h_j(x^\nu)|\}\}\\[5pt]
& & - \max_{i,j}\{(g_i(x^\nu) +
 \nabla g_i(x^\nu)\trt d(x^\nu))_+, |h_j(x^\nu) + \nabla h_j(x^\nu)\trt d(x^\nu)|\}\\[5pt]
& \overset{(c)}{\le} & \displaystyle \max_{i,j}\{(g_i(x^\nu) - g_i(x^\nu) - \nabla g_i(x^\nu)\trt d(x^\nu))_+, |h_j(x^\nu) - h_j(x^\nu) - \nabla h_j(x^\nu)\trt d(x^\nu)|\}\\[5pt]
& \le & \left\|\begin{pmatrix}
\nabla g(x^\nu)\trt\\
\nabla h(x^\nu)\trt	
\end{pmatrix}
 d(x^\nu)\right\|_\infty \le \left\|\begin{pmatrix}
\nabla g(x^\nu)\trt\\
\nabla h(x^\nu)\trt	
\end{pmatrix}\right\|_\infty \|d(x^\nu)\|,
\end{array}
\end{equation*}
where (a) follows from  $\tilde g(d(x^\nu); x^\nu) \le \kappa(x^\nu) e^m$, $-\kappa(x^\nu) e^p \le \tilde h(d(x^\nu); x^\nu) \le \kappa(x^\nu) e^p$, and $\max\{0, \alpha_1\} \le \max\{0, \alpha_2\}$ for any $\alpha_1, \, \alpha_2 \in \mathbb R$ such that $\alpha_1 \le \alpha_2$; (b) is due to A5, A7 and A9; and (c) follows from $\max\{0, \alpha_1, \beta_1\} - \max\{0, \alpha_2, \beta_2\} \le \max\{0, \alpha_1 - \alpha_2, \beta_1 - \beta_2\} $, for any $\alpha_1, \, \alpha_2, \, \beta_1, \, \beta_2 \in \mathbb R$, and $|\beta_1| - |\beta_2| \le |\beta_1 - \beta_2|$ for any $\beta_1, \, \beta_2 \in \mathbb R$. 
\hfill  \endproof
Leveraging Lemma \ref{th:prellem}, we can establish a key continuity property for the solution mapping $d(\bullet)$ of subproblem \eqref{eq:p_k}. Preliminarily, for the reader's convenience, we report the MFCQ for \eqref{eq:p_k} at $d \in \widetilde {\mathcal{X}}(x)$.
\[
\begin{array}{rl}
\{0\} = \Big\{(\xi, \pi) \, | \, & \xi \in N_{\mathbb R_-^m}(\tilde g(d; x) - \kappa(x) e^m), \pi \in N_{\kappa(x) \mathbb B_\infty^p}(\tilde h(d; x)),\\
& 0 \in \nabla_1 \tilde g(d;x)\xi + \nabla h(x)\pi + N_{\beta \mathbb B_\infty^n \cap (K - x)}(d)\Big\}
\end{array}
\] 

\begin{proposition}\label{th:keyemfcq}
Under Assumption A, let the eMFCQ hold at $\hat x \in K$. Then,

(i) the MFCQ holds at every point of $\widetilde {\mathcal{X}}(\hat x)$;

(ii)  a neighborhood $\mathcal V$ of $\hat x$ exists such that, for every point $x \in K \cap \mathcal V$, the function $d(\bullet)$ is continuous relative to $K$.
\end{proposition}
\proof{} (i) 
Since eMFCQ holds at $\hat x$, case (ii) in Lemma \ref{th:prellem} cannot occur. On the other hand, as for both cases
(i) and (iii) in Lemma \ref{th:prellem}, Slater's constraint qualification holds for $\widetilde {\mathcal{X}}(\hat x)$ and, since
$\widetilde {\mathcal{X}}(\hat x)$ is convex, this proves (i), see e.g. again \cite[Exercise 6.39]{RockWets98}, but now applied to subproblem (P$_{\hat x}$).

(ii) Since $\kappa(\bullet)$ is continuous  by Proposition \ref{eq:kappa lip}, thanks to A6 and (ii) in Lemma \ref{th:prelresfeas}, 
we can assert the outer semicontinuity, relative to $K$, at $\hat x$ for the set-valued mapping $\widetilde
{\mathcal{X}}(\bullet) = [\beta \mathbb B^n_\infty \cap (K - \bullet)] \cap \left\{d \in \mathbb R^n \, : \, \tilde g(d; \bullet) \le \kappa(\bullet)
e^m, \, -\kappa(\bullet) e^p \le \tilde h(d; \bullet) \le \kappa(\bullet) e^p\right\}$ (see \cite[Theorem 3.1.1]{banknon})

As for the inner semicontinuity property, we distinguish two cases. If $\max_{i,j}\{g_i(\hat x)_+, \, |h_j(\hat x)|\} > 0$, $\widetilde {\mathcal{X}}(\bullet)$, by
virtue of Lemma \ref{th:prellem} (i), A5, A6 and (ii) in Lemma \ref{th:prelresfeas}, is
also inner semicontinuous (see \cite[Theorem 3.1.6]{banknon}) at $\hat x$ relative to $K$. Else, if $\max_{i,j}\{g_i(\hat x)_+, \, |h_j(\hat x)|\} = 0$, in order to prove the inner semicontinuity of $\widetilde {\mathcal{X}}(\bullet) = \Big\{d \in \mathbb R^n \, : \, \tilde g(d; \bullet) \le \kappa(\bullet) e^m\Big\} \cap \Big\{d \in \mathbb R^n \, : \, -\kappa(\bullet) e^p \le \tilde h(d; \bullet) \le \kappa(\bullet) e^p, \, \|d\|_\infty \le \beta, \, d \in K - \bullet \Big\}$ at $\hat x$, relative to $K$, suffice it to show, thanks again to \cite[Theorem 3.1.6]{banknon}, that the mapping $\widetilde {\mathcal H}(\bullet) \triangleq \Big\{d \in \mathbb R^n \, : \, \tilde h(d; \bullet) = 0, \, \|d\|_\infty \le \beta, \, d \in K - \bullet \Big\}$ is inner semicontinuous at $\hat x$, relative to $K$. Suppose by contradiction that the latter mapping is not inner semicontinuous at $\hat x$ relative to $K$. Thus, there exists $\hat d \in \widetilde {\mathcal{H}}(\hat x)$, a sequence $x^\nu \underset{K}{\to} \hat x$ and $\varepsilon > 0$, such that, eventually,
\begin{equation}\label{eq:contrinner}
\|d^\nu - \hat d\| \ge \varepsilon, \qquad \forall \; d^\nu \in \widetilde {\mathcal{H}}(x^\nu).	
\end{equation}
Taking into account Lemma \ref{th:prellem} (iii), we know that $d \in \mathrm{int} (\beta \mathbb B^n_\infty) \cap \mathrm{rel \, int}(K - \hat x)$ exists such that $\tilde h(d; \hat x) = h(\hat x) + \nabla h(\hat x)\trt d = 0$. We also observe that the vectors $\nabla h_j(\hat x),$ $j = 1, \ldots, p$ are linearly independent since $\{0\} = \Big\{\pi \, | \, (\nabla h(\hat x) \pi)\trt w = 0, \forall w \in T_{\beta \mathbb B_\infty^n \cap (K - \hat x)}(v) \Big\}$, for every $v \in \widetilde {\mathcal{X}}(\hat x)$. As a consequence, denoting by $(\nabla h(\hat x)\trt)^\dagger \triangleq \nabla h(\hat x) \left[\nabla h(\hat x)\trt \nabla h(\hat x)\right]^{-1}$ the Moore-Penrose pseudoinverse of $\nabla h(\hat x)\trt$, on the one hand
\begin{equation}\label{eq:pseudo1}
\exists \, \hat z \in \mathbb R^n \, : \, \hat d = -(\nabla h(\hat x)\trt)^\dagger h(\hat x) + (I - (\nabla h(\hat x)\trt)^\dagger \nabla h(\hat x)) \hat z \in \beta \mathbb B^n_\infty \cap (K - \hat x),
\end{equation}
on the other hand,
\begin{equation}\label{eq:pseudo2}
\exists \, z \in \mathbb R^n \, : \, d = -(\nabla h(\hat x)\trt)^\dagger h(\hat x) + (I - (\nabla h(\hat x)\trt)^\dagger \nabla h(\hat x)) z \in \mathrm{int} (\beta \mathbb B^n_\infty) \cap \mathrm{rel \, int}(K - \hat x).
\end{equation}
Consider the direction $(1 - \tau^\nu) \hat d + \tau^\nu d$ with $(0,1) \ni \tau^\nu \downarrow 0$: by \cite[Theorem 6.1]{rockafellar1970convex}, we have for every $\nu$
\[
-(\nabla h(\hat x)\trt)^\dagger h(\hat x) + (I - (\nabla h(\hat x)\trt)^\dagger \nabla h(\hat x)) [\hat z + \tau^\nu (z - \hat z)] \in \mathrm{int} (\beta \mathbb B^n_\infty) \cap \mathrm{rel \, int}(K - \hat x)   
\]
and, eventually, for continuity reasons,
{\small{\[
d^\nu \triangleq -(\nabla h(x^\nu)\trt)^\dagger h(x^\nu) + (I - (\nabla h(x^\nu)\trt)^\dagger \nabla h(x^\nu)) [\hat z + \tau^\nu (z - \hat z)] \in \mathrm{int} (\beta \mathbb B^n_\infty) \cap \mathrm{rel \, int}(K - x^\nu),   
\]}}     
with $h(x^\nu) + \nabla h(x^\nu)\trt d^\nu = 0$, since the vectors $\nabla h_j(x^\nu)$, $j=1, \dots, p$ are linearly independent. Thus, $d^\nu \in \widetilde{\mathcal{H}}(x^\nu)$ exists such that $d^\nu \to \hat d$, in contradiction to \eqref{eq:contrinner}. In turn, thanks to A1, the continuity (relative to $K$) of $d(\bullet)$, leveraging \cite[Theorem 4.3.3]{banknon}, follows from \cite[Corollary 5.20]{RockWets98}.
\hfill  \endproof
To enforce the convergence results in the next section, we need
$d(\bullet)$ to be not only continuous, but also H\"older continuous on compact sets: for this reason, we introduce Assumption B.

\medskip \noindent {\bf{Assumption B}}

\medskip

\noindent {\textit{For any compact set $S\subseteq K$, two positive constants $\mu$ and $\alpha$ exist such that}}
$$
\| d(y) - d(z)\| \leq \mu \|y-z\|^\alpha, \quad \forall y, z \in S.
$$
Since it is not immediately obvious when this condition is satisfied, below we give a set of simple sufficient conditions on
$\tilde f$ and $\tilde g$ for Assumption B to hold.

\noindent {\bf{Assumption C}}
{\begin{description}
\item [{\textit{C1)}}] {\textit{$\tilde f(\bullet;\bullet)$ is locally Lipschitz continuous on $O_d \times O_x$;}}
\item [{\textit{C2)}}] {\textit{each $\tilde{g}_j(\bullet;\bullet)$ is locally Lipschitz continuous on $O_d \times O_x$.}}
\end{description}}
Note that Assumption C is automatically satisfied if we use the quadratic/linear approximations \eqref{eq:quad_lin approx}. The following proposition shows the desired result.
\begin{proposition}\label{th:uspropofupgpre}
Under Assumptions A and C, let $S\subseteq K$ be compact. Suppose further that the eMFCQ
holds
at every $\hat x \in K$. Then, there exists $\mu > 0$ such that, for every ${y},{z}\in S$,
\begin{equation}\label{eq:hol}
\|d({y})-d({z})\|\le\mu\|{y}-{z}\|^{\frac{1}{2}}.
\end{equation}
\end{proposition}
\proof{}
Preliminarily, observe  that by Proposition \ref{eq:kappa lip}, $\kappa(\bullet)$ is locally Lipschitz continuous. Furthermore, by Proposition \ref{th:keyemfcq} (i), the MFCQ holds at every point in $\widetilde{\cal X}(\hat x)$. In turn, by \cite[Theorem 3.2]{Rock85}, for every $\hat x \in K$, the set-valued mapping $\widetilde{\mathcal X}$ has the Aubin property relative to $K$ at $\hat x$ for any element belonging to $\widetilde{\mathcal X}(\hat x)$ (see \cite{RockWets98} for the definition of the Aubin property). Even more, in the light of \cite[Theorems 9.38 and 9.30]{RockWets98} being $\widetilde {\mathcal X}$ outer semicontinuous and locally bounded at $\hat x$ relative to $K$, for every $\hat x \in K$, $\widetilde{\mathcal X}$ is Lipschitz continuous (see \cite{RockWets98} for the definition of the Lipschitz continuity in the context of set-valued mappings) on a neighborhood of $\hat x$ relative to $K$. Therefore, in view of \cite[Theorem 3.3]{li2014holder}, for every $\hat x \in K$, there exist $\hat \mu > 0$ and a neighborhood ${\mathcal{V}}$ of $\hat x$ such that, for every $y, \, z \in {\mathcal{V}}\cap K$
\begin{equation*}
\|d({y})-d({z})\|\le \hat \mu \|{y}-{z}\|^{\frac{1}{2}}.
\end{equation*}
This together with the compactness of set $S$ implying a uniform bound $\mu$ across $\hat{x}$ implies \eqref{eq:hol}.
\hfill  \endproof
We recall that the KKT conditions for problem \eqref{eq:p_k} can be written as follows:
\begin{equation*}\label{eq:kktindicator}
0 \in \nabla_1 \tilde f(d(x); x) + \partial q(x + d(x)) + \nabla_1 \tilde g(d(x); x)  \xi + \nabla h(x) \pi  + N_{\beta \mathbb B^n_\infty \cap (K - x)}(d(x)),
\end{equation*}
with the KKT multipliers $\xi$ and $\pi$ satisfying the conditions $\xi \in N_{\mathbb R^m_-}(\tilde g(d(x); x) - \kappa(x) e^m)$, and $\pi \in N_{\kappa(x) \mathbb B^p_\infty} (\tilde h(d(x); x))$, respectively.

\begin{proposition}\label{lem:boundedness}
Under Assumption A, let $\hat x\in K$ and suppose that $\hat d \in \mathrm{int} (\beta \mathbb B^n_\infty) \cap \mathrm{rel \, int}(K - \hat x)$ exists such that $\tilde g(\hat d; \hat x) < \kappa(\hat x) e^m$ and either $\kappa(\hat x) e^p < \tilde h(\hat d; \hat x) < \kappa(\hat x) e^p$ (if $\kappa(\hat x) > 0$), or $\tilde h(\hat d; \hat x) = 0$ with $\{0\} = \Big\{\pi \, | \, (\nabla_1 \tilde h(\hat d; \hat x) \pi)\trt w = 0, \forall w \in T_{\beta \mathbb B_\infty^n \cap (K - \hat x)}(v) \Big\}$, for every $v \in \widetilde {\mathcal{X}}(\hat x)$ (if $\kappa(\hat x) = 0$).
Then, a neighborhood $\mathcal V$ of $\hat x$ exists such that, for every point $x \in K \cap \mathcal V$, the unique solution $d(x)$ of \eqref{eq:p_k} is a KKT point for problem \eqref{eq:p_k} and the set-valued mapping of the KKT multipliers is locally bounded at $\hat x$ relative to $K$.
\end{proposition}

\proof{}
The condition $\tilde g(\hat d; \hat x) < \kappa(\hat x) e^m$ and either $\kappa(\hat x) e^p < \tilde h(\hat d; \hat x) < \kappa(\hat x) e^p$ (if $\kappa(\hat x) > 0$), or $\tilde h(\hat d; \hat x) = 0$ and $\{0\} = \Big\{\pi \, | \, (\nabla_1 \tilde h(\hat d; \hat x) \pi)\trt w = 0, \forall w \in T_{\beta \mathbb B_\infty^n \cap (K - \hat x)}(v) \Big\}$, for every $v \in \widetilde {\mathcal{X}}(\hat x)$, (if $\kappa(\hat x) = 0$), with $\hat d \in \mathrm{int} (\beta \mathbb B^n_\infty) \cap \mathrm{rel \, int}(K - \hat x)$ is nothing else but the Slater's CQ for problem (P$_{\hat x}$), which obviously implies that the MFCQ holds at the unique solution of problem (P$_{\hat x}$) (see also again point (i) in Proposition \ref{th:keyemfcq}). The derivation of the result is then rather classical and follows from, e.g., \cite[Proposition 5.4.3]{FacchPangBk} taking into account Lemma \ref{th:prelresfeas} (ii),  Propositions \ref{eq:kappa lip} and \ref{th:keyemfcq}, and the outer semicontinuity of $N_{\beta \mathbb B_\infty^n \cap (K - \bullet)}(\bullet)$ and $N_{\kappa(\bullet) \mathbb B_\infty^p}(\bullet)$, see Lemma \ref{th:prelresfeas} (iii) and (iv) (with $\psi \equiv \kappa$), respectively. 
\hfill  \endproof

 \section{Convergence of the Method in the General Case}\label{Sec:convergence}
We are now ready to study convergence of Algorithm \ref{algoBasic}. The main result is stated below.  
\begin{theorem}\label{th:converg}
Consider the sequence $\{x^\nu\}$ generated by Algorithm \ref{algoBasic} with $\tilde f$ and $\tilde g$ satisfying Assumption A. 
The whole sequence $\{x^\nu\}$ is contained in $K$ and either is unbounded or satisfies the following assertions:
\begin{description}
\item[\rm (i)]
at least one limit limit point $\hat x$ of $\{x^\nu\}$ is generalized stationary for problem \eqref{eq:startpro}; in particular, if the eMFCQ holds at $\hat x$, then $\hat x$ is a KKT point for problem \eqref{eq:startpro};
\item[\rm (ii)] if, in addition, the eMFCQ holds at every limit point of $\{x^\nu\}$, under Assumption B, every limit point of  $\{x^\nu\}$ is a KKT solution for problem \eqref{eq:startpro}.
\end{description}
\end{theorem}
\proof{}
Since the starting point $x^0$ belongs to the convex set $K$, the stepsize satisfies $\gamma^\nu \leq 1$ by construction and, by the last constraint in (P$_{x^\nu}$), $x^\nu + d(x^\nu)\in K$ for all $\nu$, it is easily seen that all points $x^\nu$ generated by the algorithm belong to $K$.
We now assume, without loss of generality, that the sequence $\{x^\nu\}$ is bounded.

Preliminarily, observe that, at each step, the solution $d(x^\nu)$ of subproblem (P$_{x^\nu}$)
is also a KKT point for (P$_{x^\nu}$). In fact, suppose that at a certain iteration $\bar \nu$, $d(x^{\bar \nu})$ does not satisfy the KKT conditions for
(P$_{x^{\bar \nu}}$). The subproblem is always feasible by construction; let us analyze the three exhaustive cases considered in Lemma \ref{th:prellem}. In case (i), Slater's condition holds for (P$_{x^{\bar \nu}}$) and $d(x^{\bar \nu})$ is a KKT point. In case (ii), $x^{\bar \nu}$ is an ES point of \eqref{eq:startpro}: hence, we would have stopped at step \ref{S.11}. In case (iii), either Slater's condition holds for (P$_{x^{\bar \nu}}$) and $d(x^{\bar \nu})$ is a KKT point, or $x^{\bar \nu}$ is a FJ point for \eqref{eq:startpro}, in which case we would have stopped at step \ref{S.11}. Therefore, $d(x^{\bar \nu})$ is a KKT point and multipliers $\{\xi^\nu\}$ and $\{\pi^\nu\}$ exist with $\xi^\nu \in N_{\mathbb R^m_-}(\tilde g(d(x^\nu); x^\nu) - \kappa(x^\nu) e)$, $\pi^\nu \in N_{\kappa(x^\nu) \mathbb B_\infty^p}(\tilde h(d(x^\nu); x^\nu))$ and
\begin{equation}\label{eq:kktnutris}
\begin{array}{l}
0 \in \nabla_1 \tilde f(d(x^\nu); x^\nu) + \partial q(x^\nu + d(x^\nu)) + \nabla_1 \tilde g(d(x^\nu); x^\nu) \xi^\nu + \nabla h(x^\nu) \pi^\nu\\[5pt]
\qquad + N_{\beta \mathbb B^n_\infty \cap (K - x^\nu)}(d(x^\nu)).
\end{array}
\end{equation}
Using   to A1 and A4, we have
\begin{equation}\label{eq:convscobjter}
\begin{array}{rcl}
\nabla_1 \tilde f(d(x^\nu); x^\nu)\trt d(x^\nu) & = & [\nabla_1 \tilde f(d(x^\nu); x^\nu) - \nabla_1 \tilde f(0; x^\nu) + \nabla_1 \tilde f(0; x^\nu)]\trt d(x^\nu)\\[0.5em]
& \ge & c \|d(x^\nu)\|^2 + \nabla f(x^\nu)\trt d(x^\nu).
\end{array}
\end{equation}
Also, by the convexity of $q$, for every $\rho^\nu \in \partial q(x^\nu + d(x^\nu))$,
\begin{equation}\label{eq:convterm}
	{\rho^\nu}\trt d(x^\nu) \ge q(x^\nu + d(x^\nu)) - q(x^\nu). 
\end{equation}
Moreover, in view of A5, for every $i = 1, \ldots, m$,
\begin{equation}\label{eq:convscIter}
\begin{array}{rcl}
- \nabla_1 \tilde g_i(d(x^\nu); x^\nu)\trt d(x^\nu) \le \tilde g_i(0; x^\nu) - \tilde g_i(d(x^\nu); x^\nu)
\end{array}
\end{equation}
and, by A7, since $\xi^\nu$ is nonnegative, in turn,
\begin{equation}\label{eq:convscIIter}
- \xi_i^\nu \nabla_1 \tilde g_i(d(x^\nu); x^\nu)\trt d(x^\nu) \le \xi_i^\nu [\tilde g_i(0; x^\nu) - \tilde g_i(d(x^\nu); x^\nu)] = \xi_i^\nu [g_i(x^\nu) - \kappa(x^\nu)],
\end{equation}
where the equality follows observing that $\xi^\nu$ belongs to $N_{\mathbb R^m_-}(\tilde g(d(x^\nu); x^\nu) - \kappa(x^\nu) e)$. Also, taking into account $\pi^\nu \in N_{\kappa(x^\nu) \mathbb B_\infty^p}(\tilde h(d(x^\nu);x^\nu))$, we have for every $j= 1, \ldots, p$, $\pi_j^\nu [h_j(x^\nu) + \nabla h_j(x^\nu)\trt d(x^\nu)] = |\pi_j^\nu| \, \kappa(x^\nu)$ and in turn
\begin{equation}\label{eq:eqconstr}
- \pi_j^\nu \nabla h_j(x^\nu)\trt d(x^\nu) = \pi_j^\nu h_j(x^\nu) - |\pi_j^\nu| \kappa(x^\nu) \le |\pi_j^\nu| \, [|h_j(x^\nu)| - \kappa(x^\nu)].
\end{equation}
Therefore, by \eqref{eq:kktnutris}, \eqref{eq:convscobjter}, \eqref{eq:convterm}, \eqref{eq:convscIIter} and \eqref{eq:eqconstr}, we have, for some $\rho^\nu \in \partial q(x^\nu + d(x^\nu))$ and $\zeta^\nu \in N_{\beta \mathbb B^n_\infty \cap (K - x^\nu)}(d(x^\nu)),$
\begin{equation*}\label{eq:convscIIIter}
\begin{array}{l}
c \|d(x^\nu)\|^2 + \nabla f(x^\nu)\trt d(x^\nu) + q(x^\nu + d(x^\nu)) - q(x^\nu)\\[0.5em]
\qquad \le  \nabla_1 \tilde f(d(x^\nu); x^\nu)\trt d(x^\nu) + {\rho^\nu}\trt d(x^\nu)\\[0.5em]
\qquad =  - {\xi^\nu}\trt \nabla_1 \tilde g(d(x^\nu); x^\nu)\trt d(x^\nu) - {\pi^\nu}\trt \nabla h(x^\nu)\trt d(x^\nu) - {\zeta^\nu}\trt d(x^\nu)\\[0.5em]
\qquad \le  {\xi^\nu}\trt [g(x^\nu) - \kappa(x^\nu) e] + \sum_{j=1}^p |\pi_j^\nu| [|h_j(x^\nu)| - \kappa(x^\nu)]\\[0.5em]
\qquad \le  {\xi^\nu}\trt [\max_{i,j}\{g_i(x^\nu)_+, |h_j(x^\nu)|\} - \kappa(x^\nu)] e\\[0.5em]
\hspace*{32pt} + \sum_{j=1}^p |\pi_j^\nu| \, [\max_{i,j} \{g_i(x^\nu)_+, |h_j(x^\nu)|\} - \kappa(x^\nu)] =  \theta(x^\nu) \|(\xi^\nu, \pi^\nu)\|_1,
\end{array}
\end{equation*}
where the second inequality is due to $0 \in \beta \mathbb B^n_\infty \cap (K - x^\nu)$.
Hence, we get
\begin{equation}\label{eq:convscfinter}
\nabla f(x^\nu)\trt d(x^\nu) \le - c \|d(x^\nu)\|^2 +  \theta(x^\nu) \, \|(\xi^\nu, \pi^\nu)\|_1 + q(x^\nu) - 
q(x^\nu + d(x^\nu)).
\end{equation}
We also notice that, since $d(x^\nu)$ is feasible for problem (P$_{x^\nu}$), by A5, A7 and A9, for all $i$,
\begin{equation}\label{eq:convscgter}
\kappa(x^\nu) \ge \tilde g_i(d(x^\nu); x^\nu) \ge \tilde g_i(0; x^\nu) + \nabla \tilde g_i(0; x^\nu)\trt d(x^
\nu) = g_i(x^\nu) + \nabla g_i(x^\nu)\trt d(x^\nu).
\end{equation}
Let us now consider  the nonsmooth (ghost) penalty function already described in the introduction
\begin{equation}\label{eq:pendef}
W(x;\varepsilon) = f(x) + q(x) + \frac{1}{\varepsilon} \max_{i,j} \{g_i(x)_+, |h_j(x)|\},
\end{equation}
with a positive penalty parameter $\varepsilon$. This function plays a key role in the subsequent convergence 
analysis although it does not appear anywhere in the algorithm itself.

In the following analysis we will freely invoke some properties of function $(\bullet)_+ \triangleq \max\{0, 
\bullet\},$ namely $\max\{0, \alpha_1\} \le \max\{0, \alpha_2\}$ for any $\alpha_1, \, \alpha_2 \in \mathbb 
R$ such that $\alpha_1 \le \alpha_2,$ $\max\{0, a \, \alpha\} = a \, \max\{0, \alpha\}$ for any $\alpha \in 
\mathbb R$ and nonnegative scalar $a,$  and $\max\{0, \alpha_1 + \alpha_2\} \le \max\{0, \alpha_1\} + 
\max\{0, \alpha_2\}$.
We have
\begin{align}
& W(x^{\nu + 1};\varepsilon) - W(x^\nu;\varepsilon) = f(x^\nu + \gamma^\nu d(x^\nu)) - f(x^\nu) + q(x^\nu + \gamma^\nu d(x^\nu)) - q(x^\nu) \nonumber\\
& \quad  + \frac{1}{\varepsilon} \displaystyle \Big[\max_{i,j}\{g_i(x^\nu+ \gamma^\nu d(x^\nu))_+, |h_j(x^\nu+\gamma^\nu d(x^\nu))|\} \displaystyle - \max_{i,j} \{g_i(x^\nu)_+,|h_j(x^\nu)|\}\Big] \nonumber\\
& \overset{(a)}{\le}  \gamma^\nu \nabla f(x^\nu)\trt d(x^\nu) + \frac{(\gamma^\nu)^2 L_{\nabla f}}{2} \|d(x^\nu)\|^2 + q(x^\nu + \gamma^\nu d(x^\nu)) - q(x^\nu) \nonumber\\
& \quad + \frac{1}{\varepsilon} \displaystyle \Big[\max_{i,j} \{(g_i(x^\nu) + \gamma^\nu \nabla g_i(x^\nu)\trt d(x^\nu))_+, |h_j(x^\nu) + \gamma^\nu \nabla h_j(x^\nu)\trt d(x^\nu)|\}\nonumber\\  
& \quad - \displaystyle \max_{i,j}\{g_i(x^\nu)_+, |h_j(x^\nu)|\} + \frac{(\gamma^\nu)^2  \max_{i,j}\{L_{\nabla g_i}, L_{\nabla h_j}\}}{2} \|d(x^\nu)\|^2\Big] \nonumber\\[5pt]
& \overset{(b)}{\le}  \gamma^\nu \nabla f(x^\nu)\trt d(x^\nu) + \frac{(\gamma^\nu)^2}{2} (L_{\nabla f} + \frac{ \max_{i,j}\{L_{\nabla g_i}, L_{\nabla h_j}\}}{\varepsilon}) \|d(x^\nu)\|^2\nonumber\\
& \quad + q(x^\nu + \gamma^\nu d(x^\nu)) - q(x^\nu)\nonumber\\
& \quad + \frac{1}{\varepsilon} \displaystyle \Big[\max_{i,j} \{(1 - \gamma^\nu) g_i(x^\nu)_+ + \gamma^\nu \kappa(x^\nu), (1 - \gamma^\nu) |h_j(x^\nu)| + \gamma^\nu \kappa(x^\nu)\} \nonumber\\
& \quad - \displaystyle \max_{i,j}\{g_i(x^\nu)_+, |h_j(x^\nu)|\}\Big] \nonumber\\
& =  \gamma^\nu \nabla f(x^\nu)\trt d(x^\nu) + \frac{(\gamma^\nu)^2}{2} (L_{\nabla f} + \frac{ \max_{i,j}\{L_{\nabla g_i}, L_{\nabla h_j}\}}{\varepsilon}) \|d(x^\nu)\|^2\nonumber\\
& \quad + q(x^\nu + \gamma^\nu d(x^\nu)) - q(x^\nu)\nonumber\\
& \quad + \frac{1}{\varepsilon} \displaystyle \Big[(1 - \gamma^\nu) \max_{i,j} \{g_i(x^\nu)_+, |h_j(x^\nu)|\} - \displaystyle \max_{i,j}\{g_i(x^\nu)_+, |h_j(x^\nu)|\} + \gamma^\nu \kappa(x^\nu) \Big] \nonumber\\
& = \gamma^\nu \nabla f(x^\nu)\trt d(x^\nu) + \frac{(\gamma^\nu)^2}{2} (L_{\nabla f} + \frac{ \max_{i,j}\{L_{\nabla g_i}, L_{\nabla h_j}\}}{\varepsilon}) \|d(x^\nu)\|^2  \nonumber \\
& \quad + q(x^\nu + \gamma^\nu d(x^\nu)) - q(x^\nu) - \frac{\gamma^\nu}{\varepsilon} \displaystyle \Big[\max_{i,j} \{g_i(x^\nu)_+, |h_j(x^\nu)|\} - \kappa(x^\nu) \Big] \nonumber\\
&  = \gamma^\nu \nabla f(x^\nu)\trt d(x^\nu) + \frac{(\gamma^\nu)^2}{2} (L_{\nabla f} + \frac{ \max_{i,j}\{L_{\nabla g_i}, L_{\nabla h_j}\}}{\varepsilon}) \|d(x^\nu)\|^2 \nonumber - \frac{\gamma^\nu}{\varepsilon} \theta(x^\nu) \nonumber\\
& \quad q(x^\nu + \gamma^\nu d(x^\nu)) - q(x^\nu)
\label{eq:nonsmmaj1ter}\end{align}
where (a) follows applying the descent lemma to $f$, $g_i$ and $h_j$ for every $i = 1, \ldots, m$, $j = 1, \ldots, p$, with $L_{\nabla f}$, $L_{\nabla g_i}$ and $L_{\nabla h_j}$ being the Lipschitz moduli of $\nabla f$, $\nabla g_i$ and $\nabla h_j$ on the bounded set containing all iterates, and noticing that, for all $j$s,
\[|h_j(x^\nu + \gamma^\nu d(x^\nu))| \le |h_j(x^\nu) + \gamma^\nu \nabla h_j(x^\nu)\trt d(x^\nu)| + \frac{(\gamma^\nu)^2 \max_{i,j}\{L_{\nabla g_i}, L_{\nabla h_j}\}}{2} \|d(x^\nu)\|^2;
\]
(b) holds for any positive $\gamma^\nu \le 1$ since, in view of \eqref{eq:convscgter}, $\nabla g_i(x^\nu)\trt d(x^\nu) \le \kappa(x^\nu) - g_i(x^\nu)$, and, recalling $|h_j(x^\nu) + \nabla h_j(x^\nu)\trt d(x^\nu)| \le \kappa(x^\nu)$, $|h_j(x^\nu) + \gamma^\nu \nabla h_j(x^\nu)\trt d(x^\nu)| \le \gamma^\nu \kappa(x^\nu) + (1 - \gamma^\nu) |h_j(x^\nu)|$.
Furthermore, we observe that
\begin{align}\label{eq:nonsmmaj2ter}
& \gamma^\nu [\nabla f(x^\nu)\trt d(x^\nu) - \frac{1}{\varepsilon} \, \theta(x^\nu)] + q(x^\nu + \gamma^\nu d(x^\nu)) - q(x^\nu) \le \gamma^\nu [- c \|d(x^\nu)\|^2 \nonumber\\
& \quad + (\|(\xi^\nu, \pi^\nu)\|_1 - \frac{1}{\varepsilon}) \, \theta(x^\nu)] +
 \gamma^\nu [q(x^\nu) - q(x^\nu + d(x^\nu))] + q(x^\nu + \gamma^\nu d(x^\nu)) - q(x^\nu) \nonumber\\
& \le \gamma^\nu [- c \|d(x^\nu)\|^2 + (\|(\xi^\nu, \pi^\nu)\|_1 - \frac{1}{\varepsilon}) \theta(x^\nu)] 
\end{align}
where the first inequality is due to \eqref{eq:convscfinter} and the second relation is a consequence of the convexity of $q$.
We also notice that, for any fixed $x^\nu$ and for any $\eta \in (0, 1]$, there exists $\bar \varepsilon^\nu > 0$ such that
\begin{equation}\label{eq:unifsuffdescfixter}
\nabla f(x^\nu)\trt d(x^\nu) - \frac{1}{\varepsilon} \, \theta(x^\nu) + \frac{1}{\gamma^\nu} [q(x^\nu + \gamma^\nu d(x^\nu)) - q(x^\nu)] \le - \eta c \|d(x^\nu)\|^2 \qquad \forall \varepsilon \in (0, \bar \varepsilon^\nu].
\end{equation}
We now distinguish two cases.

(I) Suppose that \eqref{eq:unifsuffdescfixter} does not hold uniformly for every $x^\nu$, that is $\eta \in (0,1]$ 
and a subsequence $\{x^\nu\}_{\mathcal N}$ exists, where $\mathcal{N}\subseteq\{0, 1,2, \ldots\}$, such that 
we can construct a corresponding subsequence $\{\varepsilon^\nu\}_{\mathcal N} \in \mathbb R_+$ with $
\varepsilon^\nu \downarrow 0$ on $\mathcal N$ and
\begin{equation}\label{eq:contrter}
\nabla f(x^\nu)\trt d(x^\nu) - \frac{1}{\varepsilon^\nu} \, \theta(x^\nu) + \frac{1}{\gamma^\nu} [q(x^\nu + 
\gamma^\nu d(x^\nu)) - q(x^\nu)] > -\eta c \|d(x^\nu)\|^2
\end{equation}
for every $\nu \in {\mathcal N}$.
For \eqref{eq:contrter} to hold, relying on \eqref{eq:nonsmmaj2ter}, the multipliers' subsequence $\{(\xi^\nu, 
\pi^\nu)\}_{\mathcal N}$ must be unbounded. Combining \eqref{eq:nonsmmaj2ter} and \eqref{eq:contrter}, we 
get
\begin{equation*}\label{eq:unb}
0 \le c (1 - \eta) \|d(x^\nu)\|^2 < \left((m+p)\|(\xi^\nu, \pi^\nu)\|_\infty - \frac{1}{\varepsilon^\nu}\right) \, 
\theta(x^\nu),
\end{equation*}
and, thus, $\theta(x^\nu) > 0$ for every $\nu \in \mathcal N$.
By the previous relation and \eqref{eq:contrter}, we also have
\begin{equation}\label{eq:unb2}
\frac{1}{\varepsilon^\nu} < \frac{\nabla f(x^\nu)\trt d(x^\nu) + \frac{1}{\gamma^\nu} [q(x^\nu + \gamma^
\nu d(x^\nu)) - q(x^\nu)] + \eta c \|d(x^\nu)\|^2}{\, \theta(x^\nu)}.
\end{equation}
As $\varepsilon^\nu \downarrow 0$ on $\mathcal N$, the right hand side of \eqref{eq:unb2} goes to infinity: 
since, by the (local) Lipschitz continuity of $q$, the numerator is bounded, we have
\begin{equation}\label{eq:mustbezero}
\, \theta(x^\nu) \underset{\mathcal N}{\to} 0.
\end{equation}
Let $\hat x$ be a cluster point of the subsequence $\{x^\nu\}_{\mathcal N}$. By \eqref{eq:mustbezero}, only 
cases (ii) and (iii) in Lemma \ref{th:prellem} can occur at $\hat x \in K$. The existence of a $d$ as stipulated in 
Lemma \ref{th:prellem} (iii) would entail, by Proposition \ref{lem:boundedness}, the boundedness of the KKT 
multipliers $(\xi^\nu, \pi^\nu)$ for $\nu \in {\mathcal N}$ large enough, thus giving a contradiction. 
Therefore, by Lemma \ref{th:prellem} (ii), we conclude that $\hat x$ is either an ES or FJ point for 
\eqref{eq:startpro}.

(II) As opposed to (I), consider the case in which relation \eqref{eq:unifsuffdescfixter} holds uniformly for every 
$x^\nu$: that is, for any $\eta \in (0, 1]$, there exists $\bar \varepsilon > 0$ such that
\begin{equation}\label{eq:unifsuffdescter}
\nabla f(x^\nu)\trt d(x^\nu) - \frac{1}{\varepsilon} \, \theta(x^\nu) + \frac{1}{\gamma^\nu} [q(x^\nu + 
\gamma^\nu d(x^\nu)) - q(x^\nu)] \le - \eta c \|d(x^\nu)\|^2 \enspace \forall \varepsilon \in (0, \bar 
\varepsilon], \;\; \forall \nu.
\end{equation}

\noindent
Combining relations \eqref{eq:nonsmmaj1ter} and \eqref{eq:unifsuffdescter}, we get
\begin{equation}\label{eq:nonsmmaj3ter}
\begin{array}{rcl}
W(x^{\nu + 1}; \tilde \varepsilon) - W(x^\nu; \tilde \varepsilon) & \le & - \gamma^\nu \eta c \|d(x^\nu)\|^2 
+ \frac{(\gamma^\nu)^2}{2} (L_{\nabla f} + \frac{\max_{i,j}\{L_{\nabla g_i}, L_{\nabla h_j}\}}{\tilde \varepsilon}) 
\|d(x^\nu)\|^2\\[5pt]
& = & -\gamma^\nu \left[\eta c -  \frac{\gamma^\nu}{2} (L_{\nabla f} + \frac{\max_{i,j}\{L_{\nabla g_i}, 
L_{\nabla h_j}\}}{\tilde \varepsilon})\right] \|d(x^\nu)\|^2,
\end{array}
\end{equation}
for any $\tilde \varepsilon \in (0, \bar \varepsilon]$.
Since $\lim_\nu \gamma^{\nu}=0$, there exists a positive constant $\omega$ such that, by 
\eqref{eq:nonsmmaj3ter}, for $\nu\ge\bar{\nu}$ sufficiently large,
\begin{equation}\label{eq:nonsmmaj2bister}
W(x^{\nu+1}; \tilde \varepsilon) -  W(x^\nu; \tilde \varepsilon) \le - \omega \gamma^\nu \|d(x^\nu)\|^2.
\end{equation}
With $W$ being bounded from below, by \eqref{eq:nonsmmaj2bister}, the sequence $\{W(x^\nu; \tilde 
\varepsilon)\}$ converges and
\begin{equation*}\label{eq:summable_seriester}
\lim_{\nu}\sum_{t=\bar{\nu}}^{\nu}\gamma^{t}\|d(x^t)\|^{2}<+\infty.
\end{equation*}
Therefore, since $\sum_{\nu=0}^{\infty}\gamma^{\nu}=+\infty$, we have
\begin{equation}\label{eq:liminf is 0}
\lim \inf_{\nu \to \infty} \|d(x^\nu)\| = 0.
\end{equation}
Recalling relation \eqref{eq:thetadelta}, taking the limit on a subsequence ${\mathcal N}$ such that $\|d(x^\nu)
\| \underset{\mathcal N}{\to} 0$, we have $\theta(x^\nu) \underset{\mathcal N}{\to} 0.$
Finally, let $\hat x$ be a cluster point of subsequence $\{x^\nu\}_{\mathcal N}$. Since $\theta(x^\nu) 
\underset{\mathcal N}{\to} 0$ implies $\kappa (\hat x) = \max_{i,j} \{g_i(\hat x)_+, |h_j(\hat x)|\}$, cases (ii) 
or (iii) in Lemma \ref{th:prellem} may occur: specifically, $\hat x$ is either an ES, or a FJ, or a KKT point for 
\eqref{eq:startpro}. In particular, if the eMFCQ holds at $\hat x$, case (ii) in Lemma \ref{th:prellem} is ruled out 
and $\max_{i,j} \{g_i(\hat x)_+, |h_j(\hat x)|\}$ cannot be strictly positive; then, $\kappa (\hat x) = \max_{i,j} 
\{g_i(\hat x)_+, |h_j(\hat x)|\} = 0$. Furthermore, taking the limit in \eqref{eq:kktnutris}, we obtain, by A3, A4, 
A6-A9, KKT multipliers' boundedness and outer semicontinuity property of $\partial q(\bullet)$, and of the 
normal cone mappings $N_{\beta \mathbb B^n_\infty \cap (K - \bullet)}(\bullet)$ and $N_{\kappa(\bullet) 
\mathbb B_\infty^p}(\bullet)$ (see Lemma \ref{th:prelresfeas} (iii) and (iv) with $\psi \equiv \kappa$) and 
$N_{\mathbb R^m_-}(\bullet)$,
\[
\begin{array}{rcl}
- \nabla f(\hat x) - \nabla g(\hat x) \hat \xi - \nabla h(\hat x) \hat \pi \in \partial q(\hat x) + N_{\beta 
\mathbb B^n_\infty \cap (K - \hat x)}(0) & = & \partial q(\hat x) + \{0\} + N_{K - \hat x} (0)\\[5pt]
& = & \partial q(\hat x) + N_K(\hat x),
\end{array}
\]
with $\hat \xi \in N_{\mathbb R^m_-}(g(\hat x) - \kappa(\hat x) e) = N_{\mathbb R^m_-}(g(\hat x))$, $\hat 
\pi \in N_{\{0\}}(h(\hat x)) = \mathbb R^p$ and where the first equality follows from Lemma 
\ref{th:prelresfeas} (i). In turn, $\hat x$ is a KKT point for problem \eqref{eq:startpro}. This concludes the proof 
of case (i).

As for point (ii), observe that if, instead of the weaker \eqref{eq:liminf is 0},
\begin{equation}\label{eq:lim is 0}
\lim_{\nu \to \infty} \|d(x^\nu)\| = 0
\end{equation}
holds, we can reason similarly to what done above after \eqref{eq:liminf is 0} for any convergent subsequence 
of $\{x^\nu\}$, and conclude that (ii) holds.
Therefore, it is enough to show that Assumption B entails \eqref{eq:lim is 0}.

Consider now the compact set containing all iterates $x^\nu$.
While $\liminf_{\nu\rightarrow\infty} \|d(x^\nu)\|=0$, suppose by contradiction that $\limsup_{\nu\rightarrow\infty} \|d(x^\nu)\|>0$. Then, there exists
$\delta>0$ such that $
 \|d(x^\nu)\|>\delta$ and $ \|d(x^\nu)\|<\delta/2
$
for infinitely many $\nu$s. Therefore, there is an infinite
subset of indices ${\cal N}$ such that, for each $\nu\in{\cal N}$,
and some $i_{\nu}>\nu$, the following relations hold:
\begin{equation}\label{eq:con1}
 \|d(x^\nu)\|<\delta/2,\hspace{6pt}\|d({x}^{i_{\nu}})\|>\delta
\end{equation}
and, if $i_\nu > \nu + 1$,
\begin{equation}
\delta/2\le\|d({x}^{j})\|\le\delta,\hspace{6pt}\nu<j<i_{\nu}.\label{eq:con2}
\end{equation}
Hence, for all $\nu\in{\cal N}$, we can write
\begin{equation}\label{eq:ineqser}
\begin{array}{rcl}
\delta/2 & < &  \|d({x}^{i_{\nu}})\|-\|d({x}^{\nu})\|
 \, \leq \, \|d({x}^{i_{\nu}}) - d({x}^{\nu})\| \overset{(a)}{\le} \mu \|{x}^{i_{\nu}}-{x}^{\nu}\|^\alpha\\[0.5em]
 & \overset{(b)}{\le} & \mu \left[\sum_{t=\nu}^{i_{\nu}-1}\gamma^{t} \|d(x^t)\|\right]^\alpha \overset{(c)}{\le} \mu \delta^\alpha\left(\sum_{t=\nu}^{i_{\nu}-1}\gamma^{t}\right)^\alpha,
\end{array}
\end{equation}
where (a) is due to Assumption B with $\alpha$ and $\mu$ positive scalars, (b) comes from the triangle inequality and the updating rule of the algorithm and in (c) we used
\eqref{eq:con2}. By \eqref{eq:ineqser} we have
\begin{equation}\label{eq:absu}
\underset{_{\nu\rightarrow\infty}}\liminf \;\; \mu \delta^\alpha\left(\sum_{t=\nu}^{i_{\nu}-1}\gamma^{t}\right)^\alpha>0.
\end{equation}
We prove next that \eqref{eq:absu} is in contradiction with the convergence of $\{W(x^\nu;\tilde \varepsilon)\}$
for any $\tilde \varepsilon \in (0, \bar \varepsilon]$, where $\bar\varepsilon$ is defined around
\eqref{eq:unifsuffdescter}. To this end, we first
show that $\|d({x}^{\nu})\|\ge\delta/4$, for sufficiently large $\nu\in{\cal N}$. Reasoning as in \eqref{eq:ineqser}, we have
\begin{equation*}
\|d({x}^{\nu+1})\| - \|d({x}^{\nu})\| \le \mu \|
{x}^{\nu+1}-{x}^{\nu}\|^\alpha \le \mu (\gamma^\nu)^\alpha \|d({x}^{\nu})\|^\alpha,
\label{eq:ineqser2}
\end{equation*}
for any given $\nu$. For  $\nu\in{\cal N}$ large enough so that $\mu ({\gamma^\nu})^\alpha(\delta/4)^\alpha<\delta/4$, suppose by contradiction that $\|
d({x}^{\nu})\|<\delta/4$; this would give $\|d({x}^{\nu+1})\|<\delta/2$ and, thus,
condition \eqref{eq:con2} (or \eqref{eq:con1}) would be violated. Then, it must be
$
\|d({x}^{\nu})\|\ge\delta/4.
$
From this, and using \eqref{eq:nonsmmaj2bister}, we have, for sufficiently large $\nu\in{\cal N}$,
\begin{equation}\label{eq:fin}
W({x}^{i_{\nu}}; \tilde \varepsilon) \le W({x}^{\nu}; \tilde \varepsilon) -\omega\sum_{t=\nu}^{i_{\nu}-1}\gamma^{t}\|d({x}^{t})\|^{2} \le W({x}^{\nu}; \tilde \varepsilon) - \omega\frac{\delta^{2}}{16}\sum_{t=\nu}^{i_{\nu}-1}\gamma^{t}.
\end{equation}
Since $\{W(x^\nu; \tilde \varepsilon)\}$ converges, as established above immediately after \eqref{eq:nonsmmaj2bister},
renumbering if necessary, relation \eqref{eq:fin} implies $\sum_{t=\nu}^{i_{\nu}-1}\gamma^{t} \to 0$, in contradiction with \eqref{eq:absu}.
This shows that \eqref{eq:lim is 0} holds and concludes the proof of the theorem. \hfill
\hfill   \endproof

\begin{remark}\label{rem:boundedness}  Convergence to a generalized stationary point is obtained in Theorem \ref{th:converg}  if the sequence $\{x^\nu\}$ is bounded.  In our framework, generating
an unbounded sequence is a natural possibility that cannot and should not  be excluded in principle, since
we do not make  any standard  assumption such as feasibility, regularity of
the constraints, coercivity or even existence of a stationary point.
Of course, the question arises of when the sequence generated by the algorithm is bounded; can we give {\em a priori} conditions that guarantee the boundedness of the iterations?
It is possible to give a satisfactory answer to this question, at the price of a much more convoluted analysis; we refer the reader to \cite{facchinei2017ghost} for developments in this direction. 
Here we only mention that if $K$ is bounded, a case very frequent in applications,  the sequence 
$\{x^\nu\}$, which is contained in $K$, is also surely bounded.
\end{remark}


\section{Convergence of a Simplified Version of the Method for Problems with Convex Constraints satisfying eMFCQ}\label{sec:newalgo}
The result in Theorem \ref{th:converg} is very broadly applicable, implying some notion of stationarity for limit points even for problems satisfying essentially no structural assumptions. In this section, we add two classical assumptions with respect to the constraints, namely, convexity, and a standard constraint qualification.

A simple corollary of Theorem \ref{th:converg} is obtained assuming that the eMFCQ holds everywhere; note that this is commonly assumed in papers considering convergence of algorithms for nonlinear constrained optimization. Theorem \ref{th:converg} immediately gives the following result.

\begin{corollary}\label{cor:converg}
Consider the sequence $\{x^\nu\}$ generated by Algorithm \ref{algoBasic} with $\tilde f$ and $\tilde g$ such that Assumption A holds and suppose that the eMFCQ holds everywhere in $K$. Then, the whole sequence $\{x^\nu\}$ is contained in $K$ and either the sequence $\{x^\nu\}$ is unbounded or the following assertions hold:
\begin{description}
\item[\rm (i)]
at least one limit limit point $\hat x$ of $\{x^\nu\}$  is a KKT point for problem \eqref{eq:startpro};
\item[\rm (ii)] if  Assumption B also holds, then every limit point of  $\{x^\nu\}$ is a KKT point for problem \eqref{eq:startpro}.
\end{description}
\end{corollary}
In what follows we further assume that the constraints are convex and we explore the consequences of this structural property: we show that Algorithm \ref{algoBasic} can actually be simplified while stronger convergence results can be obtained.
Therefore, from now on we make the following assumption.
\medskip 

\noindent {\bf{Assumption D}}

\medskip

\noindent {\textit{Each $g_i$ is convex and $h$ is linear, i.e. $h(x) = Ax + b$, for some $p \times n$ matrix $A$ and vector $b \in \mathbb R^p$, on an open neighborhood of $K$, and $K$ is bounded. Furthermore, the feasible set ${\mathcal{X}}$ of \eqref{eq:startpro} is nonempty and the eMFCQ holds at every point in $K$.
}}
\medskip

\noindent
Note that $f$ is not assumed to be convex. In this setting, we show that one can always set $\kappa(x) =0$, thus avoiding the non negligible task of computing this quantity, and make the resulting modified version of Algorithm \ref{algoBasic} resemble a pure, classical SQP-type method.
The key point here is that $\kappa(x)$ is introduced to relax the constraints in subproblem \eqref{eq:p_k} in order to ensure nonemptiness of its feasible set and some continuity properties of $d(x)$. It turns out that, alternatively, Assumption D is also sufficient to achieve these results. The search direction we consider is now defined to be the solution of the following strongly convex optimization problem: 
\begin{equation}\label{eq:pc_k}
\begin{array}{cl}
\underset{d}{\mbox{minimize}} &  \tilde f(d;x) + q(x + d)\\
\mbox{s.t.} & g(x) + \nabla g(x)\trt d  \le 0\\[5pt]
&  A (x+d) + b = 0,\\[5pt]
& d \in K - x,
\end{array}\tag{P$_{x}^c$}
\end{equation}
whose feasible set is denoted by $\widetilde {\mathcal X}^c(x)$.
Note that this is just a particular case of \eqref{eq:p_k} where we take $\tilde g$ to be the linear approximations of $g$ and omit the constraint $\|d\|_\infty \le \beta$; this latter simplification is possible because we assume
$K$ to be bounded and therefore, if we take $\beta$ to be larger than the diameter of $K$, this constraint is superfluous (i.e. never active) and can be omitted.
Note also that if $x^0$ satisfies the linear equality constraints, every point generated by the algorithm will satisfy them, and the equality constraints in \eqref{eq:pc_k} can be rewritten as $Ad =0$, since at each iteration $b+ Ax^\nu =0$. 

Our first order of business is then to ensure that the feasible set of \eqref{eq:pc_k} is always nonempty; this is rather classical to show, even if in our setting not totally immediate.

\begin{proposition}\label{pro:nonempty}
Under Assumption D, for any $x\in K$  
there exists $d \in \mathrm{rel \, int}(K - x)$ such that $g(x) + \nabla g(x)\trt d < 0$ and $A(x + d) + b = 0$. Furthermore $\{0\} = \Big\{\pi \, | \, (A\trt \pi)\trt w = 0, \forall w \in T_{K - x}(v) \Big\}$ for every $v \in \widetilde{\mathcal X}^c(x)$. A fortiori, the feasible set of \eqref{eq:pc_k} is nonempty.
\end{proposition}

\proof{}
The proof is an adaptation of the one for Lemma \ref{th:prellem} (iii). Take a feasible point $\hat x$. In view of convexity, we recall that the eMFCQ holds at $\hat x$ if and only if (see \cite[Exercise 6.39 (b)]{RockWets98}) conditions \eqref{eq:EMFCQ1}, which under Assumption D read as follows, hold: 
\begin{gather}
	\{0\} = \Big\{\pi \, | \, (A\trt \pi)\trt w = 0, \forall w \in T_K(\hat x) \Big\}, \nonumber \\[7pt]
 \exists \hat d \,\in\, \mathrm{rel \, int} \, T_K(\hat x):\; \nabla g_i(\hat x)\trt \hat d < 0, \quad \forall i: g_i(\hat x) = 0, \qquad A \hat d = 0. \label{eq:rw}
\end{gather}
First, observing that the tangent cone to the convex set $K$ at a point $\hat x$ is given by the closure of the cone of feasible directions for $K$ at $\hat x$, and borrowing again from the proof of \cite[Theorem 6.9]{RockWets98}, we have $\mathrm{rel \, int} \, T_K(\hat x) = \left\{d \in \mathbb R^n \, | \, \exists \, \alpha > 0 \enspace \text{with} \enspace \hat x + \alpha d \in \mathrm{rel \, int} K \right\}$ due to \cite[Proposition 2.40]{RockWets98}. Hence, in view of \cite[Theorem 6.1]{rockafellar1970convex}, for every $\tau > 0$ sufficiently small, $\hat x + \tau \hat d \triangleq \tilde x \in \mathrm{rel \, int} K$ as well.
As a consequence, by \eqref{eq:rw},  $\tilde x$ is still feasible for \eqref{eq:startpro} and, because of 
$\nabla g_i(\hat x)\trt \hat d < 0, \; \forall i$ such that $g_i(\hat x) = 0$, and simple continuity arguments, it holds that
$g(\tilde x)$ is stricty feasible, i.e., $g(\tilde x) <0$.

Concerning the feasible set $\widetilde {\mathcal{X}}^c(x)$ of subproblem \eqref{eq:pc_k}, take $d \triangleq \tilde x - x$. Clearly, $d \in \mathrm{rel \, int}(K - x)$. Furthermore it can readily be seen that $\tilde x$ satisfies the linear constraints, i.e., $b + A(x + (\tilde x - x)) = 0$ holds. Moreover, by using convexity, we can write, for all $i$,
\[
0 > g_i(\tilde x) \geq g_i(x)  + \nabla g_i(x)\trt(\tilde x - x) =  g_i(x)  + \nabla g_i(x)\trt d.
\]
Finally, following the same line of reasoning (by contradiction) in the proof of Lemma \ref{th:prellem} (iii) (see in particular the developments below \eqref{eq:linearindenh}), but here considering the constraint $\|d\|_\infty\le \beta$ as never active and thus omitted, since the eMFCQ holds everywhere, we get, for every $v \in \widetilde{\mathcal X}^c(x)$,
\begin{equation}\label{eq:linearindenhsimpl}
\{0\} = \Big\{\pi \, | \, (A\trt \pi)\trt w = 0, \forall w \in T_{K - x}(v) \Big\}.	
\end{equation} \hfill   \endproof
Paralleling Proposition \eqref{th:keyemfcq} (ii), we now establish that $d(x)$, also as the unique solution of the modified subproblem \eqref{eq:pc_k}, remains continuous. Since the feasible set of \eqref{eq:pc_k}  is nonempty and we assume, as usual, that the objective function of that subproblem is strongly convex (see Assumption A), $d(x)$ is well-defined to begin with.

\begin{proposition}\label{pro:continuity modified}
Under Assumptions A and D, the function $d(\bullet)$ is continuous on $K$.
\end{proposition}
\proof{}
The proof can be derived from the one for Proposition \ref{th:keyemfcq}. More specifically, in view of the continuity of the functions involved and by (ii) in Lemma \ref{th:prelresfeas}, the set-valued mapping $\widetilde
{\mathcal{X}}^c(\bullet) = (K - \bullet) \cap \left\{d \in \mathbb R^n \, : \, g(\bullet) + \nabla g(\bullet)\trt d \le 0, A(\bullet + d) + b = 0\right\}$ is outer semicontinuous relative to $K$ at any $x \in K$, thanks to \cite[Theorem 3.1.1]{banknon}. To show that $\widetilde {\mathcal{X}}^c(\bullet) = \Big\{d \in \mathbb R^n \, : \, g(\bullet) + \nabla g(\bullet)\trt d \le 0\Big\} \cap \Big\{d \in \mathbb R^n \, : \, A (\bullet + d) + b = 0, \, d \in K - \bullet \Big\}$ at $\hat x$ is also inner semicontinuous at any $x \in K$, relative to $K$, it suffices to prove, thanks again to \cite[Theorem 3.1.6]{banknon} and Proposition \ref{pro:nonempty}, that the mapping $\widetilde {\mathcal H}^c(\bullet) \triangleq \Big\{d \in \mathbb R^n \, : \, A (\bullet + d) + b = 0, \, d \in K - \bullet \Big\}$ is inner semicontinuous at $x$, relative to $K$. In the light of Proposition \ref{pro:nonempty}, this can be done by a  reductio ad absurdum, following the same steps as in the proof of Proposition \ref{th:keyemfcq} (ii) (see in particular the developments below relation \eqref{eq:contrinner}), by recalling that here the constraint $\|d\|_\infty\le \beta$ is never active and thus omitted.

Finally, thanks to A1, the continuity (relative to $K$) of $d(\bullet)$, leveraging \cite[Theorem 4.3.3]{banknon}, follows from \cite[Corollary 5.20]{RockWets98}.
 \hfill \endproof
 
 The following proposition shows that, classically, $\|d(x)\|$ is a stationarity measure.
 \begin{proposition}\label{pro:d is 0}
Under Assumptions A and D, $d(x)=0$ if and only if $x$ is a KKT point of \eqref{eq:startpro}.
\end{proposition}
\proof{}
 The (unique) solution of the strongly convex subproblem \eqref{eq:pc_k}  satisfies the KKT conditions because of Proposition \ref{pro:nonempty}. This KKT conditions can be written as
 \begin{equation}\label{eq:KKTsub}
\begin{array}{c}
0 \in \nabla_1 \tilde f(d;x) + \partial q(x+d) + \nabla g(x) \xi + A\trt \pi + N_{K-x}(d) \\[5pt]
0\leq \xi \perp \big(g(x) + \nabla g(x)\trt d\big)\leq 0\\[5pt]
A(x+d) + b =0\\[5pt]
d\in K -x,
\end{array}
\end{equation}
where, we recall, the variable is $d$. Taking into account A4 and A9, and the fact that 
$N_{K-x}(0) = N_K(x)$, the assertion of the proposition can be checked easily by comparing 
the KKT system \eqref{eq:KKTstartpro} for the original problem to that for the subproblem, i.e. \eqref{eq:KKTsub}.
 \hfill \endproof
 We now have all the required preliminary derivations to analyze the following algorithm, which is  a variant of Algorithm \ref{algoBasic} where we use \eqref{eq:pc_k} instead of \eqref{eq:p_k} to compute the direction and therefore avoid the calculation of $\kappa$.
\IncMargin{1em}
\begin{algorithm}
\KwData{$\gamma^\nu \in (0,1]$ such that \eqref{eq:gamma} holds, $x^{0} \in K$,  $\nu \longleftarrow 0$\;}
\Repeat{
{\nlset{(S.1)} \If{$x^{\nu}$ {\emph{is a KKT solution for}} \eqref{eq:startpro}}{{\bf{stop}} and {\bf{return}} $x^\nu$\;} \label{S.11c}}
\nlset{(S.2)} compute the solution $d(x^\nu)$ of problem (P$^c_{x^\nu}$)\; \label{S.12c}
\nlset{(S.3)} set $x^{\nu+1}=x^{\nu}+\gamma^{\nu}d(x^\nu)$, $\nu\longleftarrow\nu+1$\; \label{S.13c}}{}
\caption{\label{algoBasicC} Simplified DSM Algorithm for \eqref{eq:startpro} under Assumption D}
\end{algorithm}
\DecMargin{1em} 
We can now establish the convergence properties of this modified algorithm. Note that in the following theorem below we assume that Assumption C holds, but this places a requirement only on $\tilde f$, since $\tilde g$ is the linearization of $g$ and therefore Assumption C2 is automatically satisfied.
\begin{theorem}\label{th:convergC}
Consider the sequence $\{x^\nu\}$ generated by Algorithm \ref{algoBasicC} under Assumptions A, C and D. The whole sequence $\{x^\nu\}$ is bounded and  contained in $K$ and  each of its limit points is a KKT solution for problem \eqref{eq:startpro}.
\end{theorem}
\proof{}
The proof is formally similar to that of Theorem \ref{th:converg} and we do not repeat it for sake of brevity.
It is enough to follow the proof of Theorem \ref{th:converg} step by step, setting $\kappa(x) =0$ wherever this quantity appears and taking into account these points, which primarily even simplify the analysis further:
\begin{itemize}
\item by Assumption D, the sequence $\{x^\nu\}$ is bounded, since $K$ is bounded;
\item Proposition \ref{pro:nonempty} shows that $d(x)$ is well-defined;
\item Proposition \ref{pro:continuity modified}  gives the required continuity of $d(x)$ which, in the proof of Theorem \ref{th:converg}, was guaranteed by the results in Sections \ref{sec:stationary} and \ref{Sec:dirfind}; 
\item at the beginning of the proof we need to show that $d(x)$ is actually a KKT point for subproblem \eqref{eq:pc_k} and this is now guaranteed by Proposition \ref{pro:nonempty}, where we prove that that Slater's CQ holds for subproblem \eqref{eq:pc_k};
\item since $\kappa(x) =0$, $\theta(x)$ reduces to nothing else but the classical measure of the violation of the constraints: $\theta(x)  =  \max_{i,j} \{g_i(\hat x)_+, |h_j(\hat x)|\}$;
\item case I in the proof of of Theorem \ref{th:converg} cannot now occur because of the eMFCQ assumption;
\item in case II there, we have $\liminf_{\nu \to \infty} \|d(x^\nu)\| =0$. Together with the continuity of $d(x)$ and Proposition  \ref{pro:d is 0}, this implies that there exists at least one limit point of the sequence generated by the algorithm that turn out to be a KKT solution for \eqref{eq:startpro};
\item Assumption C1 guarantees that actually  every limit point of the bounded sequence $\{x^\nu\}$ is a KKT solution of \eqref{eq:startpro} since it implies Assumption B.
\end{itemize}
\hfill \endproof

\section*{Funding}
Francisco Facchinei was partially supported by Progetto di Ateneo Distributed optimization algorithms for Big Data.
Vyacheslav Kungurtsev was supported by the OP VVV project CZ.02.1.01/0.0/0.0/16\_019/0000765 ``Research
Center for Informatics".
Lorenzo Lampariello was partially supported by the MIUR PRIN 2017 (grant 20177WC4KE).
 Gesualdo Scutari was partially supported by the NSF   Grants CIF 1564044, CIF 1719205, and CMMI 1832688; and the ARO under Grant W911NF1810238.

\bibliographystyle{abbrv}
\bibliography{Surbib}

\end{document}